\documentclass[aos,MSNbibl,seceqn,dvips]{arximspdf}

\doi{10.1214/14-AOS1215} 
\volume{42}
\issue{3}
\pubyear{2014}
\firstpage{1102}
\lastpage{1130}

\makeatletter
\newcommand{\eqref}[1]{(\ref{#1})}
\renewcommand{\citep}[1]{(\citeauthor{#1} \citeyear{#1})}
\def\m{\mathcal}
\def\mb{\mathbb}
\def\mr{\mathrm}
\def\mrr{\operatorname}

\def\tr{\operatorname{tr}}
\newtheorem{thmm}{Theorem}[section]
\newtheorem{lem}[thmm]{Lemma}

\newproclaim{ass}[thmm]{Assumption}

\newproclaim{defn}[thmm]{Definition}
\newproclaim{rem}[thmm]{Remark}

\newproclaim{remark}{Remark}
\def\T{\mathrm{T}}
\newcommand{\Real}{\mathbb{R}}

\newcommand\bbP{\mathbb{P}}
\newcommand\bbE{\mathbb{E}}

\makeatother

\begin{document}
\begin{frontmatter}

\title{Posterior contraction in sparse Bayesian factor models for
massive covariance matrices}
\runtitle{Bayesian Covariance matrix estimation}

\begin{aug}
\author[a]{\fnms{Debdeep}~\snm{Pati}\corref{}\ead[label=e1]{debdeep@stat.fsu.edu}\thanksref{t1}},
\author[b]{\fnms{Anirban}~\snm{Bhattacharya}\ead[label=e2]{anirbanb@stat.tamu.edu}\thanksref{t1}},
\author[c]{\fnms{Natesh~S.}~\snm{Pillai}\ead[label=e3]{pillai@fas.harvard.edu}\thanksref{t2}}
\and
\author[d]{\fnms{David} \snm{Dunson}\ead[label=e4]{dunson@stat.duke.edu}\thanksref{t1}}
\runauthor{Pati, Bhattacharya, Pillai and Dunson}
\thankstext{t1}{Supported in part by the DARPA MSEE program and
Grant number R01 ES017240-01 from the
National Institute of Environmental Health Sciences (NIEHS) of the National
Institutes of Health (NIH).}
\thankstext{t2}{Supported in part by NSF Grant DMS-11-07070.}

\affiliation{Florida State University, Texas A\&M University,
Harvard University\break and Duke University}

\address[a]{D. Pati\\
Department of Statistics\\
Florida State University\\
Tallahassee, Florida\\
USA\\
\printead{e1}}
\address[b]{A. Bhattacharya\\
Department of Statistics\\
Texas A\&M University\\
College Station, Texas\\
USA\\
\printead{e2}\hspace*{20.5pt}}
\address[c]{N.~S. Pillai\\
Department of Statistics\\
Harvard University\\
Cambridge, Massachusetts\\
USA\\
\printead{e3}}
\address[d]{D. Dunson\\
Department of Statistical Science\\
Duke University\\
Durham, North Carolina\\
USA\\
\printead{e4}}
\end{aug}

\received{\smonth{9} \syear{2013}}
\revised{\smonth{2} \syear{2014}}

%
\begin{abstract}
Sparse Bayesian factor models are routinely implemented for
parsimonious dependence modeling and dimensionality reduction in
high-dimensional applications. We provide theoretical understanding of
such Bayesian procedures in terms of posterior convergence rates in
inferring high-dimensional covariance matrices where the dimension can
be larger than the sample size.
Under relevant sparsity assumptions on the true covariance matrix, we
show that commonly-used point mass mixture priors on the factor
loadings lead to consistent estimation in the operator norm even when
$p \gg n$. One of our major contributions is to develop a new class of
continuous shrinkage priors and provide insights into their
concentration around sparse vectors. Using such priors for the factor
loadings, we obtain similar rate of convergence as obtained with point
mass mixture priors.
To obtain the convergence rates, we construct test functions to
separate points in the space of high-dimensional covariance matrices
using insights from random matrix theory; the tools developed may be of
independent interest. We also derive minimax rates and show that the
Bayesian posterior rates of convergence coincide with the minimax rates
upto a $\sqrt{\log n}$ term.

%
\end{abstract}

%

\end{frontmatter}

\section{Introduction}\label{sec1}

It is now routine to collect data where the dimension $p$ is much
larger than the sample size $n$, and interest focuses on the covariance
structure. In this context, even a simple parametric model like the
Gaussian distribution leads to a high-dimensional model space and it
becomes necessary to reduce the effective number of parameters via
imposing sparsity or some lower-dimensional structure. Sparse Bayesian
factor models \cite{West03} provide one popular choice in
applications, but currently lack theoretical support. In this paper, we
close this gap by studying asymptotic properties for scenarios in which
$p$ grows faster than $n$.

Factor models \cite{bartholomew1987latent} aim to explain dependence
among multivariate observations through shared dependence on a smaller
number of latent factors. Given $n$ i.i.d. observations $y_i \in
\mathbb
{R}^p$, a latent factor model is given by
%
\begin{equation}
\label{eq:modeldef_intro} y_i = \Lambda\eta_i +
\varepsilon_i,\qquad  \varepsilon_i \sim
\mathrm{N}_p(0, \Omega), i = 1,\ldots,n,
\end{equation}
where $\Lambda$ is a $p \times k$ factor loadings matrix with $k \ll
p$, $\eta_i \sim\mathrm{N}_k(0,\mathrm{I})$ are standard normal latent
factors, and $\varepsilon_i$ is a residual having diagonal covariance
$\Omega= \operatorname{diag}(\sigma_1^2,\ldots,\sigma_p^2)$.
Marginalizing out the latent factors, $y_i \sim\mathrm{N}_p(0,\Sigma)$ with
%
\begin{equation}
\label{eq:factor_decomp} \Sigma= \Lambda\Lambda^{\T} + \Omega,
\end{equation}
so that the right-hand side has at most $p(k+1)$ parameters compared to
$O(p^2)$ parameters in an unstructured covariance matrix.

A prior distribution on $(\Lambda, \Omega)$ induces a prior
distribution on $\Sigma$ and we are interested in studying
concentration of the corresponding posterior measure around a ``true''
covariance matrix in \emph{operator norm} when the dimensionality $p =
p_n$ can be much larger than the sample size $n$. This setting has
motivated abundant frequentist work, with rates of convergence of
various regularized covariance estimators derived in \cite
{bickel2008regularized,Bickel08,lam2009sparsistency,el2008operator,cai2011adaptivethr,bunea2012sample}
among others. Minimax optimal rates for specific sparsity classes have
also been derived in \cite{cai2010optimal,cai2012optimal}.
There is a relatively smaller but increasing literature on asymptotic
properties of Bayesian procedures in models with growing dimension,
primarily focused on linear or generalized linear models; refer to
\cite
{ghosal1999as,ghosal2000as,belitser2003adaptive,armagan2011generalized,bontemps2011bernstein,castilloneedles}
among others. To the best of our knowledge, the present paper is the
first to study the
asymptotic properties of Bayesian covariance estimation via factor
models in the $p_n \gg n$ regime.

We now summarize the main results obtained in this paper. Although the
original specification of the factor model reduces the number of
parameters from quadratic to linear in $p_n$, the estimation problem is
still challenging when $p_n \gg n$. To address this challenge, \cite
{West03} introduced \emph{sparse factor modeling} to allow many of the
loadings to be exactly equal to zero through a point mass mixture prior
having a probability mass at zero; see also \cite{Lucasetal06,Carvalho08} for modifications and applications in genomics. Recently,
\cite{castilloneedles} studied posterior concentration in estimating a
sparse high-dimensional mean using such point mass mixture priors.
However, it is not clear whether the \emph{induced} prior on the
covariance from such sparsity favoring priors on the factor loadings
would lead to consistent covariance estimation in the $p_n \gg n$
setting. We answer the question in the affirmative and derive the rate
of convergence of the posterior in Section~\ref{sec:main_res},
explicitly characterizing the dependence on the dimensionality $p_n$,
the true number of factor $k_{0n}$, the column sparsity $s_n$ in the
true loadings and the growth rate of the largest eigenvalue $c_n$ of
the true covariance. In particular, the dimensionality enters the rate
through a logarithmic factor, providing justification of usage of such
methods in \emph{ultra high-dimensional} settings. It may be remarked
here that the usual practice of assuming the eigenvalues of the true
covariance to be bounded is restrictive in our context and we relax
that assumption.

Although point mass mixture priors are amenable to incorporate
sparsity, exploring the model space via MCMC can be daunting and may
lead to slow mixing and convergence of the algorithm \cite
{polson2010shrink}. To address such problems through block updating,
while allowing a weaker notion of sparsity in which elements are close
to zero instead of exactly zero, continuous shrinkage priors can be
used. Such priors have become common in regression \cite
{park2008bayesian,carvalho2010horseshoe,hans2011elastic,armagan2011generalized},
with \cite{polson2010shrink} providing a unifying local--global scale
mixture representation. Although computationally attractive, the lack
of tight concentration bounds for such priors has limited the study of
their asymptotic properties. One of our main contributions is to
develop a novel class of continuous shrinkage priors and derive
nonasymptotic bounds on the concentration and dimensionality of such
priors. Based on these results, we show that the proposed continuous
shrinkage prior leads to the same rate of posterior convergence as the
point mass mixture priors in estimating large covariance
matrices.\looseness=1

The Birg{\'e}--Le Cam testing theory \cite
{birge1984tests,le1986asymptotic} for the Hellinger metric is commonly
used in Bayesian asymptotics \cite{ghosal2000convergence} to separate
points in the parameter space. However, generalization of the testing
argument to other norms has been relatively unexplored. A notable
exception is \cite{gine2011rates} who advocated the use of
concentration inequalities based on empirical process techniques to
derive tests in the $L_r$ metric in a nonparametric function estimation
context. See also \cite{ray2013inverse} for an usage of concentration
bounds for centered linear estimators in the context of test
construction in Bayesian inverse problems. In the setting of large
covariance estimation in \emph{operator norm}, we construct tests inspired
by results from the nonasymptotic theory of random matrices, which
might be of independent interest in related settings.

Finally, we use Fano's lemma to derive the \emph{minimax rate} of
convergence for the class of covariance matrices considered in this
paper and show that the posterior indeed convergences at the minimax
rate up to a $\sqrt{\log n}$ term.

There is a sizeable literature studying asymptotic properties of
various aspects of factor analysis, including consistent estimation of
factor loadings and latent factors \cite{bai2003inferential} and the number of factors \cite
{bai2002determining,lam2011factor}. Fan, Fan and Lv [19] studied rates
of convergence of high-dimensional covariance estimates based on factor
models, with~\cite{fan2011high} extending their results to approximate
factor models that allow nondiagonal $\Omega$ in~\eqref
{eq:factor_decomp}. This work assumes that the factor scores $\eta_i$
are known, while we consider the fundamentally different setting in
which the factor scores are unknown while also studying concentration
of a Bayesian posterior instead of convergence of a point estimate.

The rest of the paper is organized as follows. After setting up the
basic notation and definitions in Section~\ref{sec:prelim}, we state
our assumptions and their implications in Section~\ref{sec:ASS}.
In Section~\ref{sec:prior}, we discuss our prior distributions. The
main results of this paper are stated in Section~\ref{sec:main_res}.
Section~\ref{sec:sims} contains some numerical simulations.
In Section~\ref{sec:disc}, we prove a number of concentration bounds
for the shrinkage prior introduced in Section~\ref{sec:prior}, while in
Section~\ref{sec:test}, we elucidate our test construction. These
results are used to prove the main results in Section~\ref{sec:main_pf}.
Proof of some technical lemmas are given in a supporting document.

\section{Preliminaries} \label{sec:prelim}

Given sequences $a_n, b_n$, we shall denote $a_n = O(b_n)$ or $a_n
\lesssim b_n$ if there exists a global constant $C$ such that $a_n \leq
C b_n$. Similarly, we define
$a_n \gtrsim b_n$ and $a_n \asymp b_n$.

Given a metric space $(X,d)$, let $N(\varepsilon; X, d)$ denote its
$\varepsilon$-covering number, that is, the minimum number of balls of
radius $\varepsilon$ needed to cover $X$.

For a vector $ x \in\mathbb{R}^r$, $\Vert x\Vert_2$ denotes its Euclidean
norm. We will use $\mathcal{S}^{r-1}$ to denote the unit Euclidean
sphere $\{x \in\mathbb{R}^r \dvtx \Vert x\Vert_2 = 1\}$ and $\Delta
^{r-1}$ to
denote the $(r-1)$-dimensional simplex $\{ x = (x_1, \ldots,
x_{r})^{\T
} \dvtx x_j \geq0, \sum_{j=1}^{r} x_j = 1\}$. Further, let $\Delta
_0^{r-1}$ denote $\{ x = (x_1, \ldots, x_{r-1})^{\T} \dvtx x_j \geq0,
\sum_{j=1}^{r-1} x_j \leq1\}$.

For a square matrix $A$, $\operatorname{tr}(A)$ and $\vert A\vert$, respectively,
denote the trace and the determinant of $A$. For a $p \times r$ matrix
$A = (a_{jj'})$ with $p \geq r$, let
$s_{(1)} \geq s_{(2)} \geq\cdots\geq s_{(r)} \geq0$ denote the
singular values of $A$ (or equivalently the eigenvalues of $\sqrt {A^{\T
}A}$) arranged in decreasing order. We shall use $s_{\min}(A)$ and
$s_{\max}(A)$ to denote the smallest and largest singular values, respectively.
The Frobenius norm ($\Vert\cdot\Vert_F$) and the operator norm
($\Vert\cdot\Vert_2$) are defined in the usual way, with $ \Vert
A\Vert_F:= \sqrt {\operatorname{tr}(A^{\T} A)}$ and $\Vert A\Vert_2:= \sup_{x \in\mathcal
{S}^{r-1}} \Vert A x\Vert_2 = s_{\max}(A)$. Also $\Vert A\Vert_1 =
\sum_{j=1}^p
\sum_{h=1}^r |A_{jh}|$ is the $l_1$ norm of $\operatorname{vec}(A)$.
We will derive posterior convergence rates in the \emph{operator norm}.

For a subset $S \subset\{1, \ldots, p\}$, let $|S|$ denote the
cardinality of $S$ and define $\theta_S = (\theta_j \dvtx j \in S)$ for a
vector $\theta\in\mathbb{R}^p$.
Denote $\operatorname{supp}(\theta)$ to be the \emph{support} of $\theta$,
that is, the subset $S_0 \subset\{1, \ldots, p\}$ corresponding to the
nonzero entries of $\theta$. We shall continue to use the same notation
for a subset of entries and support for matrices~$\Lambda$, where it
has to be interpreted that $\Lambda$ is vectorized column-wise. Let
$l_0[s; p]$ be the space of $s$-sparse vectors $\theta\in\mathbb
{R}^p$ with\break  $|\operatorname{supp}(\theta)| \leq s$.

Throughout $C, C'$ are generically used to denote positive constants
whose values might change from one line to the next but are independent
from everything else.

Finally, let $\mathcal{C}_n$ denote the cone of covariance matrices of
size $p_n \times p_n$ and let $\Sigma_{0n} \in\mathcal{C}_n$ denote a
true sequence of covariance matrices.\setcounter{footnote}{2}\footnote{As a convention, we make
the dependence of all quantities on $n$ explicit, and only omit that in
a few places for notational convenience.} We observe
\[
y_1, \ldots, y_n \stackrel{\mathrm{i.i.d}} {\sim}
\mathrm{N}_{p_n}(0, \Sigma_{0n})
\]
and set $\mathbf{y}^{(n)} = (y_1, \ldots, y_n)$. We model the data as
%
\begin{equation}
\label{eq:factor_model} y_i \stackrel{\mathrm{i.i.d}} {\sim}
\mathrm{N}_{p_n}(0, \Sigma_n),\qquad \Sigma_n =
\Lambda_n \Lambda_n^{\T} + \Omega_n,
\Omega_n = \sigma ^2 \mathrm{I}_{p_n}.
\end{equation}

We will denote our prior distribution on $\mathcal{C}_n$ (constructed
in Section~\ref{sec:prior}) by $\Pi_n(\cdot)$ and the corresponding posterior
distribution by $\Pi_n(\cdot| \mathbf{y}^{(n)})$. 
\section{Assumptions}\label{sec:ASS}
In this section, we state our assumptions on the true data generating
model and briefly discuss their implications.
Let $\Real^{p \times k}$ denote the class of real-valued $p \times k$
matrices. We start with the following assumptions on the true
covariance matrix of the observed data $\mathbf{y}^{(n)}$.
%
\begin{ass} \label{ass:truth}
The true sequence of covariance matrices $\Sigma_{0n}$ are of the form
\[
\hspace*{-3pt}\mathrm{(A0)}\quad\Sigma_{0n} = \Lambda_{0n}
\Lambda_{0n}^{\T} + \Omega_{0n},\qquad
\Lambda_{0n} \in\Real^{p_n \times k_{0n}}, k_{0n} \leq
p_n, \Omega_{0n} = \sigma_{0n}^2
\mathrm{I}_{p_n}. 
\]
%
\end{ass}

Assumption (\ref{ass:truth}) says that the true sequence of covariances
$\Sigma_{0n}$ admit a factor decomposition as in (\ref
{eq:factor_decomp}) with $\Omega_{0n} = \sigma_{0n}^2 \mathrm
{I}_{p_n}$. We make the following assumptions on $\Lambda_{0n}$ and
$\sigma_{0n}^2$.
%
\begin{ass}\label{ass:truecov1} There exist sequences of positive real
numbers $c_n, s_n$ with $c_n \lesssim s_n$, such that:
%
\begin{longlist}[(A1)]
\item[(A1)] $\lim_{n \rightarrow\infty} c_n k_{0n}^{3/2} \sqrt {\frac
{s_n \log p_n}{n}} \sqrt{\log n} = 0; k_{0n}^{3/2} \sqrt{\frac
{s_n \log p_n}{n}} (\log n)^{3/2}= O(1)$.

\item[(A2)] Each column of $\Lambda_{0n}$ belongs to $l_0[s_n; p_n]$.

\item[(A3)] $\Vert\frac{1}{c_n} \Lambda_{0n}^{\T} \Lambda_{0n} -
\mathrm{I}_{k_{0n}}\Vert_2 = o(k_{0n} \sqrt{\log k_{0n}/n})$.

\item[(A4)] There exists a constant $\sigma_0^{(1)}$ such that
$\sigma
_0^{(1)} \leq\sigma_{0n}^2 \leq c_n$.
\end{longlist}
\end{ass}

We now discuss implications of each of the above assumptions.
\begin{itemize}
\item If $ k_{0n} = O(1)$, $c_n, s_n \asymp\log p_n$, the first part
of {(A1)} allows $p_n$ to grow faster than $n$ under the mild
assumption of
$(\log p_n)^{5}\log n /n \to0$. In this case, $p_n$ can be of the
order of $\exp(n^{\alpha})$ for any $\alpha\in(0, 1/5)$. The second
part is a very mild requirement given the first part; indeed if $p_n =
\exp(n^{\alpha})$, $s_n = n^{\beta}$ and $k_{0n} = n^{\gamma}$ for
appropriate $\alpha, \beta, \gamma> 0$ such that the first part of
(A1) holds, then the second part follows from the first.

\item In gene-expression studies, we expect each factor is related to
only a relatively
small number of variables, representing a sparse, parsimonious
structure underlying the associations among genes. Following the
motivation in \cite{West03}, usually a small number of latent factors
associate with the response so that only those genes with nonzero
loadings on those factors are relevant. This is reflected through {
(A2)}, requiring the loadings columns to be sparse with $s_n \ll p_n$
many signals per column.

\item Conditions similar to {(A3)} appear in the econometric factor
model setting \cite{fan2008high,fan2011high} referred to as
``pervasive.'' We provide an intuition based on random matrix theory
which suggests that {(A3)} is indeed mild and expected to be
satisfied by a large class of loadings.
As our emphasis is on sparse factor models, a realistic generative
model for the true loadings would be
\[
{\lambda}_{0jh} \sim(1 - \pi_n) \delta_0 +
\pi_n \mathrm{N}(0, 1),
\]
where $\lambda_{0jh} = [\Lambda_{0n}]_{jh}$, $\delta_0$ denotes a point
mass at zero and we set $\pi_n = s_n/p_n$ to reflect the sparsity
assumption in {(A2)}.
Using a modification of Theorem~5.39 of~\cite{vershynin2010introduction},
$\Vert\frac{1}{p_n} \Lambda_{0n}^{\T} \Lambda_{0n} - \pi_n
\mathrm {I}_{k_{0n}} \Vert_2 \leq C \frac{ \sqrt{k_{0n}} }{ \sqrt {p_n} } \Vert\pi _n \mathrm{I}_{k_{0n}}\Vert_2$,
or equivalently,
$\Vert\frac{1}{s_n} \Lambda_{0n}^{\T} \Lambda_{0n} - \mathrm
{I}_{k_{0n}} \Vert_2 \leq C \frac{ \sqrt{k_{0n}} }{ \sqrt{p_n} }$,
with probability at least $1 - e^{-C' k_{0n}}$. We can thus choose $c_n
= s_n$ and $\sqrt{k_{0n}/p_n}$ is smaller than $k_{0n} \sqrt{\log
k_{0n}/n}$ if $p_n = \exp(n^{\alpha})$.

\item{(A4)} simply posits an upper and lower bound on the residual
variance. The lower bound is used to avoid $\Sigma_{0n}$ being
ill-conditioned,\footnote{The constant lower bound on $\sigma_{0n}^2$
can be relaxed as long as $s_{\max}(\Sigma_{0n})/s_{\min}(\Sigma_{0n})
\lesssim n$.} while the upper bound ensures that the larger
contribution to $\Vert\Sigma_{0n}\Vert_2$ comes from the loadings
$\Lambda
_{0n}$. In particular, {(A3)} and {(A4)} imply $\| \Sigma_{0n}\|
_2 \asymp c_n$, allowing the largest eigenvalue to grow with increasing
dimension.
\end{itemize}

We denote by $\mathcal{C}_{0n}$ the class of covariance matrices
satisfying {(A0)}--{(A4)} in Assumptions \ref{ass:truth} and
\ref{ass:truecov1}. Clearly, any $\Sigma_{0n} \in\m C_{0n}$ can be
parameterized by $(k_{0n}, \Lambda_{0n}, \sigma_{0n}^2)$, where
$\Lambda
_{0n} \in\mb R^{p_n \times k_{0n}}$.

%
\section{Prior distribution} \label{sec:prior}
We consider model \eqref{eq:factor_model} with $\Lambda_n \in\Real
^{p_n \times k}$. We specify priors on the residual variance $\sigma
^2$, the number of factors $k$ and the factor loadings (conditional on
the number of factors) $\Lambda_n \mid k$ below.

For the residual variance $\sigma^2$, we assign a gamma prior
$f_\sigma
$ on $(0, \infty)$,
{\renewcommand{\theequation}{PR}
\begin{equation}
\label{eq:prior_res} \sigma^2 \sim\operatorname{Ga}(a, b).\vspace*{-9pt}
\end{equation}}

For the number of factors $k$, we assume a prior distribution $\pi_k$
which decays exponentially,
\renewcommand{\theequation}{\arabic{section}.\arabic{equation}}
\setcounter{equation}{0}
\begin{equation}
\label{eqn:priork} \pi_k(k > j) \leq\exp(-C j),
\end{equation}
for all $j \geq j_0$ for some $j_0 \in\mathbb{N}$. Additionally, assume
%
\begin{equation}
\label{eqn:priork2} \pi_k(k = k_{0n}) \geq\exp(-C
s_n k_{0n} \log n),
\end{equation}
where $s_n$ is the sequence appearing in {(A2)}. For instance, a
Poisson distribution on $k$ with rate parameter $1$ will satisfy \eqref
{eqn:priork} and \eqref{eqn:priork2} if $\log k_{0n} \leq s_n \log n$,
which is automatically satisfied given {(A1)}.

Conditional on $k$, we consider two classes of prior distributions on
the factor loadings $\Lambda_{n}$.
We first consider a class of point mass mixture priors on the loadings
similar to that advocated by \cite{West03},
\renewcommand{\theequation}{PL1}
\begin{eqnarray}\qquad
\label{eq:prior_loadings1}  \lambda_{jh} \mid k,\qquad \pi&\sim&(1 - \pi)
\delta_0 + \pi g(\cdot),\qquad j = 1, \ldots, p_n; h = 1,
\ldots, k,
\nonumber
\\[-8pt]
\\[-8pt]
\nonumber
 k &\sim&\pi_{k}, \qquad\pi\sim\operatorname{Beta}(1, \kappa k_{0n}
p_n + 1), \qquad \kappa> 0, 
\end{eqnarray}
where $\delta_0$ denotes a point mass at zero and $g$ is an absolutely
continuous density on $\mathbb{R}$ with exponential tails or heavier.

For linear models, \cite{scott2010bayes} showed that such point mass
mixture priors with a beta hyper-prior on the mixture probability lead
to an automatic multiplicity correction. \cite{jiang2007bayesian}
proved optimality results in estimating the predictive distribution
under such priors in generalized linear models accommodating diverging
numbers of predictors. Castillo and van der Vaart \cite
{castilloneedles} studied concentration properties of a class of prior
distributions similar to (\ref{eq:prior_loadings1}) on a
high-dimensional normal mean and showed that they lead to the minimax
optimal rate of convergence.

As mentioned in the \hyperref[sec1]{Introduction}, although point mass mixture priors
are conceptually appealing in allowing exact sparsity and often leading
to appealing theoretical properties, posterior computation under such
priors can be daunting in high-dimensional cases. As an alternative, a
rich variety of continuous shrinkage priors have been developed that
admit a scale mixture representation \cite{polson2010shrink}. A
fundamental hurdle in studying theoretical properties of such priors is
the difficulty of obtaining tight bounds on their concentration and
implied dimensionality. With the motivation of developing a continuous
shrinkage prior that can be shown to concentrate near sparse vectors
and approximate point mass mixture priors, we propose a novel class of
priors. We use such priors for the factor loadings, but they should be
broadly applicable in other high-dimensional settings.

Let $\operatorname{DE}(\psi)$ denote the Laplace or double-exponential density
with scale parameter $\psi$ with a density given by
\renewcommand{\theequation}{\arabic{section}.\arabic{equation}}
\setcounter{equation}{2}
\begin{equation}
\label{eqn:DEfun}
f(x) = \frac{1}{2 \psi} e^{- {\vert x\vert }/{ \psi} },\qquad x \in \mathbb{R}.
\end{equation}

Draw the elements of a high-dimensional vector $\theta\in\mathbb
{R}^p$ through the following hierarchical mechanism:
{\renewcommand{\theequation}{PS}
\begin{equation}
\label{eq:shprior_prop} \theta_j \sim\operatorname{DE}(\tau\gamma_j),\qquad
\tau\sim f_{\tau}, \gamma= (\gamma_1, \ldots,
\gamma_p)' \sim f_{\gamma}. 
\end{equation}}
\hspace*{-3pt}In \eqref{eq:shprior_prop}, $\tau> 0$ is a global scale parameter and
$\gamma\in\Delta^{p-1}$ is a vector of local scale parameters. We set
$f_{\tau}$ to be an $\operatorname{exp}(1/2)$ density. We draw $\tilde
{\gamma}
= (\gamma_1, \ldots, \gamma_{p-1})^{\T} \in\Delta_0^{p-1}$ from a
$\operatorname{Dir}(\alpha/p, \ldots, \alpha/p)$ density and set $\gamma_p
= 1
- \sum_{j=1}^{p-1} \gamma_j$. For a detailed discussion on the
properties of the prior \eqref{eq:shprior_prop}, refer to Section~\ref
{sec:disc}.

Given $k$, we consider the prior \eqref{eq:shprior_prop} on the
vectorized loadings $\operatorname{vec}(\Lambda_n) \in\mb R^{p_n k}$ as an
alternative to \eqref{eq:prior_loadings1}; note that the Dirichlet
concentration parameter becomes $\alpha/(p_n k)$ in this case.


\section{Main results}\label{sec:main_res}

With the prior specification complete, we now state the main results of
this paper. The proofs are available in Section~\ref{sec:main_pf}.
Theorems \ref{thmm:op_norm_rates} and \ref{thmm:op_norm_rates_shrinkage}
assume the true number of factors to be bounded, which is generalized
in Theorem~\ref{thmm:k_0}. Recall the class of ``true'' covariance
matrices $\m C_{0n}$ from Section~\ref{sec:ASS}. We first establish the rate of
posterior convergence in operator norm using the point mass priors
\eqref{eq:prior_loadings1} on the loadings in Theorem~\ref{thmm:op_norm_rates}.

\begin{thmm}\label{thmm:op_norm_rates}
Suppose $\Sigma_{0n} \in\mathcal{C}_{0n}$ with $s_n \gtrsim\log p_n$
and $k_{0n} = O(1)$, and model \eqref{eq:factor_model} is fitted with a
prior distribution on the number of factors satisfying~\eqref
{eqn:priork} and \eqref{eqn:priork2}. Assume independent priors $\Pi
(\Lambda\mid k )$ and $\Pi(\sigma^2)$ on the loadings and the residual
variances as in \eqref{eq:prior_loadings1} and \eqref{eq:prior_res},
respectively. Then, with $\varepsilon_n = c_n \sqrt{\frac{s_n \log
p_n}{n}} \sqrt{\log n}$ and for some constant $M > 0$,
%
\begin{equation}
\label{eq:post_prob_lim} \lim_{n \rightarrow\infty} \bbE_{\Sigma_{0n}}
\Pi_n\bigl(\Vert\Sigma _n - \Sigma_{0n}
\Vert_2 > M \varepsilon_n \mid\mathbf{y}^{(n)}
\bigr) = 0,
\end{equation}
where $\bbE_{\Sigma_{0n}}$ denotes an expectation with respect to the
joint distribution of~$\mathbf{y}^{(n)}$.
\end{thmm}
We next show in Theorem~\ref{thmm:op_norm_rates_shrinkage} that our
proposed shrinkage prior on the loadings achieves the same posterior
rate of convergence as for the point mass mixture priors.

\begin{thmm}\label{thmm:op_norm_rates_shrinkage}
Assume the same setup as in Theorem~\ref{thmm:op_norm_rates}, with the
point mass prior \eqref{eq:prior_loadings1} on the loadings replaced by
the shrinkage prior \eqref{eq:shprior_prop} on the vectorized loadings
given $k$. Then \eqref{eq:post_prob_lim} is satisfied with $\varepsilon_n
= c_n \sqrt{\frac{s_n \log p_n}{n}} \sqrt{\log n}$.
\end{thmm}
We show in Section~\ref{sec:minimax} that $c_n \sqrt{s_n \log p_n/n}$
is the minimax rate of estimating $\Sigma_{0n} \in\mathcal{C}_{0n}$ in
operator norm with $k_{0n} = O(1)$. Thus, the posterior rate of
convergence obtained in both Theorems \ref{thmm:op_norm_rates} and \ref
{thmm:op_norm_rates_shrinkage} is equal to the minimax rate up to a
$\sqrt{\log n}$ term. For a general $k_{0n}$, we establish analogous
versions of Theorems \ref{thmm:op_norm_rates} and \ref
{thmm:op_norm_rates_shrinkage} below.
%
\begin{thmm}\label{thmm:k_0}
If $\Sigma_{0n} \in\mathcal{C}_{0n}$ with $s_n k_{0n} \gtrsim\log
p_n$, the convergence rates in both Theorems \ref{thmm:op_norm_rates}
and \ref{thmm:op_norm_rates_shrinkage} are modified to $c_nk_{0n}^{3/2}
\sqrt{\frac{s_n \log p_n}{n}} \sqrt{\log n}$.\vadjust{\goodbreak}
\end{thmm}
Clearly, Theorem~\ref{thmm:k_0} permits consistent estimation in
operator norm even if $p_n = \exp(n^{\alpha})$, $s_n = n^{\beta}$ and
$k_{0n} = n^{\gamma}$ for appropriate $\alpha, \beta, \gamma> 0$. At
this point, we do not know whether the rate obtained in Theorem~\ref
{thmm:k_0} is minimax-optimal and substantial further work seems
necessary to prove such a result.

%
\subsection{A lower bound to the minimax rate} \label{sec:minimax}

Minimax optimal rates in operator norm for high-dimensional covariance
matrix estimation have been established for a class of bandable
matrices by \cite{cai2010optimal} and a class of covariance matrices
with sparse columns by \cite{cai2012optimal}. Although $\mathcal
{C}_{0n}$ has a nonempty intersection with the class $\mathcal
{G}_0(c_{n,p})$ in \cite{cai2012optimal}, there exists a large subclass
of matrices which lie in one and not in the other. Moreover, the
existing minimax results on large covariance estimation assume the
eigenvalues of the true sequence of covariance matrices to be bounded.
For example, \cite{cai2010optimal} and \cite{cai2012optimal} assume
that $y_i$ is sub-Gaussian, that is, for all $t > 0$ and $v \in\mathbb
{R}^{p}$ with $\| v \|_2 = 1$,
$\bbP(| v^{\T} (y_1 - \bbE y_1) | > t) \leq\exp(- t^2 / 2 \tau^2)$.
The parameter $\tau$ is assumed to be a constant and its role in the
rate is not characterized. For $y_1 \sim\mathrm{N}_{p}(0, \Sigma)$, a
standard tail bound for the normal distribution implies
\[
\bbP\bigl(\bigl| v^{\T} y_1 \bigr| > t\bigr) \leq\exp \biggl( -
\frac{t^2}{ 2 v^{\T}
\Sigma
v} \biggr) \leq\exp \biggl( - \frac{t^2}{ 2 \| \Sigma\|_2} \biggr).
\]
For $\Sigma= \Lambda\Lambda^{\T} + \sigma^2 \mathrm{I}_p \in
\mathcal
{C}_{0n}, \| \Sigma\|_2 = \Vert\Lambda\Vert_2^2 + \sigma^2 \asymp
c_n$ by
Assumption~\ref{ass:truecov1}, so that $\tau\asymp\sqrt{c_n}$ in our
case. Hence, the growth rate of $\| \Sigma\|_2$ needs to be accounted
for in our calculations. With this motivation, we study minimax lower
bounds for $\mathcal{C}_{0n}$ in Theorem~\ref{thmm:lb} below.
%
\begin{thmm}\label{thmm:lb}
If $\hat{\Sigma}_n$ is a sequence of estimators of $\Sigma_{0n} \in
\mathcal{C}_{0n}$ with $k_{0n} = O(1)$, then
%
\begin{equation}
\inf_{\hat{\Sigma}_n}\sup_{\Sigma_{0n} \in\mathcal{C}_{0n}} \Vert\hat {
\Sigma}_n - \Sigma_{0n} \Vert_2 \geq
c_n \sqrt{s_n \frac{\log p_n}{n}}.
\end{equation}
\end{thmm}
\begin{pf}
We will use Fano's lemma to derive a lower bound for the minimax risk.
Let $\m F$ be a parameter space of covariance matrices and we observe
$y_1, \ldots, y_n \sim\mathrm{N}(0, \Sigma)$ with $\Sigma\in\m F$. Let
$\Theta=  \{\Sigma_{(1)}, \ldots, \Sigma_{(m_n)}  \}, m_n
\geq
2$ be a finite subset of $\m F$ and let $\bbP^{(j)}$ denote the joint
distribution of $y_1, \ldots, y_n$ independently distributed as $\mathrm
{N}(0, \Sigma_{(j)})$, $1 \leq j \leq m_n$. Let $\hat{\Sigma}$ be an
estimator for $\Sigma$.
Suppose for all $j \neq j'$, we have that
\[
d(\Sigma_{(j)}, \Sigma_{(j')}) \geq d_{m_n},\qquad \mrr{KL}
\bigl(\bbP^{(j)}, \bbP^{(j')}\bigr) \leq K_{m_n}.
\]
%
Letting $\bbE_j$ denote the expectation under $\bbP^{(j)}$, Fano's
lemma (as in \cite{yu1997assouad}) implies
%
\begin{equation}
\label{eq:fan} \max_{1\leq j \leq m_n} \mathbb{E}_j d(\hat{
\Sigma}, \Sigma) \geq \frac{d_{m_n}}{2} \biggl(1 - \frac{K_{m_n} + \log2}{\log m_n} \biggr).
\end{equation}
We first introduce notation and then proceed to construct our finite
parameter set~$\Theta$. Let $q_n = p_n-1$. Define $\m M:= \{ x \in
\mb
R^{q_n} \dvtx x_j \in\{0, 1\}\ \forall  j, \|x \|_1 = s_n \}$ to be
the collection of all binary vectors of length $q_n$ with exactly $s_n$
ones. Let $\| \cdot\|_H$ denote the Hamming distance between two
binary strings, so that $\| x - y \|_H = \sum_{j=1}^{q_n} 1(x_j \neq
y_j)$. Let $b_j = (\theta_j, 0)$ denote the $p_n$-dimensional vector
obtained by appending zero at the end of $\theta_j$. With this
notation, set
\[
\Sigma_{(j)} = \beta\mr I_{p_n} + \gamma b_j
b_j^{\T} + \kappa e_{p_n} e_{p_n}^{\T},\qquad
j \in\m M,
\]
where $e_p \in\mb R^p$ is the vector with $1$ in the $p${th}
coordinate and zero elsewhere, and $\gamma< \beta< \kappa$ are
sequences to be chosen below.

We now state Lemmas \ref{lem:fan1} and \ref{lem:pack}; for clarity in
notation, we drop the subscript $n$ in both lemmata. Refer to the
\hyperref[app]{Appendix} for a proof.
%
\begin{lem}\label{lem:fan1}
For $j \neq j'$ and $0 \leq r \leq s$, if $\| \theta_j - \theta_{j'}
\|
_H = 2(s-r)$, then
\[
\Vert\Sigma_{(j)} - \Sigma_{(j')}\Vert_2 =\gamma
\sqrt{s^2 - r^2},\qquad \mrr{KL}\bigl(\bbP^{(j)},
\bbP^{(j')}\bigr) = \frac{n}{2} \frac
{t^2}{ts+1}
\bigl(s^2 - r^2\bigr),
\]
where $t = \gamma/\beta$.
\end{lem}
%
\begin{lem}\label{lem:pack}
Given $s \geq6$, there exists a subset $\m M_0 = \{\theta_1, \ldots,
\theta_m\}$ of $\m M$ with $m \asymp\exp(C s\log p)$ and $\Vert
\theta _j - \theta_{j'}\Vert_H \geq s/3$ for all $1 \leq j \neq j'
\leq m$, where
$C$ is a positive constant independent of $p$.
\end{lem}

We set $\Theta= \{ \Sigma_{(j)}\dvtx \theta_j \in\m M_0\}$. Since
$\Vert\theta_j - \theta_{j'}\Vert_H \geq s_n/3$ by Lemma~\ref
{lem:pack}, the
quantity $r_n = r_n(j, j')$ appearing in Lemma~\ref{lem:fan1} is
bounded above by $5s_n/6$ for all pairs $j \neq j' \in\m M_0$. Hence,
we can choose $d_{m_n} = C_1 \gamma s_n$ and $K_{m_n} = n (ts_n)^2 = n
(\gamma s_n/\beta)^2$ in \eqref{eq:fan}. To obtain $d_{m_n}$ as a lower
bound to the minimax risk up to a constant, we need to set
$K_{m_n}/\log m_n = C'$ for some constant $C' \in(0, 1)$. Since $\log
m_n \asymp C s_n \log p_n$, we obtain, by choosing $\beta= c_n$, that
$d_{m_n}^2 = C (\gamma s_n)^2 = C \frac{c_n^2 s_n \log p_n }{n}$ for
some absolute constant $C$.
\end{pf}

\section{Simulation studies} \label{sec:sims}
In this section, we consider a number of simulation cases to compare
our proposed continuous shrinkage prior \eqref{eq:shprior_prop} with
existing methods including the point mass priors \eqref
{eq:prior_loadings1} on the loadings matrix and the sample covariance
matrix $S = (n-1)^{-1} \sum_{i=1}^n (y_i - \bar{y}) (y_i - \bar
{y})^{\T
}$. For prior \eqref{eq:prior_loadings1}, we use a standard Laplace
distribution on the signal coefficients.

We also compare our methods with Principal Orthogonal complement
Thresholding (POET) of \cite{fan2011large} which is based on an
additive decomposition of the covariance matrix in terms of a low rank
matrix and a sparse residual covariance matrix. POET estimates the
factors and the loadings by thresholding the principal components of
the sample covariance matrix. Finally, we provide results for the
adaptive thresholding method (AT) of \cite{cai2011adaptivethr} which
thresholds the entries of the sample covariance matrix, with the
resulting thresholded estimator $\hat{\Sigma}$ being of the form
$\hat
{\Sigma}_{jj'} = S_{j j'} 1(|S_{j j'} | > \delta\kappa_{j j'}) $,
where $\delta$ is a tuning parameter and $\kappa_{j j'}$ is a threshold
specific to the corresponding entry of $S$. We chose the tuning
parameter $\delta$ by 5-fold cross-validation as suggested by \cite
{cai2011adaptivethr}. We also implemented the same procedure with the
default choice of $\delta= 2$; the results were worse in all cases,
and hence are not reported.\looseness=-1

We describe the two simulation settings below:
\begin{longlist}[1.]
\item[1.] $y_i, i = 1, \ldots, n$ are generated from $\mathrm
{N}_{p_n}( 0, \Sigma_{0n})$, where $\Sigma_{0n} = \Lambda_{0n}
\Lambda
_{0n}^{\T} + \sigma_0^2 \mathrm{I}_{p_n}$ and $\Lambda_{0n}$ is a $p_n
\times k_{0n}$ matrix with $s_n = \log p_n$ nonzero entries per column
and $k_{0n} = 1 \mbox{ or } \log p_n$. The nonzero entries were drawn
uniformly between $1$ and~$2$. These simulations were designed to mimic
assumptions (A0)--(A4) in Section~\ref{sec:ASS}.

\item[2.] This setting is designed to illustrate the performance
of our method under model misspecification. We let $\Sigma_{0n} =
\Lambda_{0n} \Lambda_{0n}^{\T} + \Omega_{0n}$, where $\Lambda
_{0n}$ is
as in simulation setting (1), but $\Omega_{0n}$ is nondiagonal,
corresponding to the covariance matrix of an autoregressive sequence
with pure error variance $0.4$ and autoregressive coefficient $0.1$.
\end{longlist}

For each simulation setting, we choose two sample sizes, namely $n =
50, 100$ and for each value of $n$, we let $p_n = 100, 200$. For each
$(n, p_n)$ pair, we consider $50$ simulation replicates. For the
Bayesian methods, the posterior mean is used as a point estimate.
Tables~\ref{tab:1} and \ref{tab:2} summarize the results across the
simulation replicates for the two simulation settings, respectively, to
compare the operator norm difference between the estimator resulting
from the different methods and the truth. In particular, the average
error across $50$ replicates is provided, with standard error in parenthesis.

\begin{table}
\tabcolsep=0pt
\caption{Simulation setting 1. Comparative performance in covariance
matrix estimation for \protect\eqref{eq:prior_loadings1}, \protect
\eqref
{eq:shprior_prop}, POET, AT. The average error in operator norm across
simulation replicates is tabulated}\label{tab:1}
{\fontsize{8}{10}\selectfont
\begin{tabular*}{\textwidth}{@{\extracolsep{\fill}}lcccccccc@{}}
\hline
\multicolumn{1}{@{}l}{\textbf{n}} &
\multicolumn{4}{c}{\textbf{50}} & \multicolumn{4}{c@{}}{\textbf{100}}
\\[-6pt]
&
\multicolumn{4}{c}{\hrulefill} & \multicolumn{4}{c@{}}{\hrulefill}\\
\multicolumn{1}{@{}l}{$\bolds{p_n}$}&
\multicolumn{2}{c}{\textbf{100}} &
\multicolumn{2}{c}{\textbf{200}} &
\multicolumn{2}{c}{\textbf{100}} &
\multicolumn{2}{c@{}}{\textbf{200}} \\[-6pt]
&
\multicolumn{2}{c}{\hrulefill} &
\multicolumn{2}{c}{\hrulefill} &
\multicolumn{2}{c}{\hrulefill} &
\multicolumn{2}{c@{}}{\hrulefill} \\
$\bolds{k_{0n}}$ & \textbf{1}&
$\bolds{\log p_n}$ & \textbf{1} &$\bolds{\log p_n}$
& \textbf{1} & $\bolds{\log p_n}$ &
\textbf{1} & $\bolds{\log p_n}$ \\
\hline
\eqref{eq:prior_loadings1} & 0.98 (0.43) & 2.84 (1.12) & 10.06 (5.68) &
\phantom{0}9.79 (4.90) & 8.83 (0.12) & 12.82 (1.42) &15.90 (0.26)& 16.07 (1.77) \\
\eqref{eq:shprior_prop} & 1.03 (0.38) & 3.95 (1.69) & \phantom{0}5.96 (1.81)
 &
\phantom{0}7.01 (2.01) & 1.74 (0.83) & \phantom{0}3.43 (1.10) &
\phantom{0}3.66 (1.83) & \phantom{0}4.21 (1.20) \\[3pt]
POET & 2.89 (0.41) & 6.98 (1.28) & \phantom{0}8.90 (2.11) &12.41 (2.69) &
 3.08 (0.64)
& \phantom{0}5.72 (1.09) &\phantom{0}7.32 (1.51) &\phantom{0}7.51 (1.44) \\
AT & 1.93 (0.57) & 4.71 (2.97) & \phantom{0}6.92 (5.43) &\phantom{0}8.86 (3.79) & 2.11 (0.71) &
\phantom{0}3.26 (1.08)& \phantom{0}3.80 (2.03) & \phantom{0}4.37 (1.34) \\
SC & 2.79 (0.36) & 7.08 (1.33) & \phantom{0}9.01 (2.22) &12.73 (2.80) & 3.06 (0.65) &
\phantom{0}5.73 (1.14) & \phantom{0}7.34 (1.52) & \phantom{0}7.52 (1.46) \\
\hline
\end{tabular*}}
%
\end{table}

\begin{table}
\tabcolsep=0pt
\caption{Simulation Setting 2. Comparative performance in covariance
matrix estimation for \protect\eqref{eq:prior_loadings1}, \protect
\eqref
{eq:shprior_prop}, POET, AT. The average error in operator norm across
simulation replicates is tabulated}\label{tab:2}
{\fontsize{8}{10}\selectfont
\begin{tabular*}{\textwidth}{@{\extracolsep{\fill}}lcccccccc@{}}
\hline
\multicolumn{1}{@{}l}{\textbf{n}} &
\multicolumn{4}{c}{\textbf{50}} & \multicolumn{4}{c@{}}{\textbf{100}}
\\[-6pt]
&
\multicolumn{4}{c}{\hrulefill} & \multicolumn{4}{c@{}}{\hrulefill}\\
\multicolumn{1}{@{}l}{$\bolds{p_n}$}&
\multicolumn{2}{c}{\textbf{100}} &
\multicolumn{2}{c}{\textbf{200}} &
\multicolumn{2}{c}{\textbf{100}} &
\multicolumn{2}{c@{}}{\textbf{200}} \\[-6pt]
&
\multicolumn{2}{c}{\hrulefill} &
\multicolumn{2}{c}{\hrulefill} &
\multicolumn{2}{c}{\hrulefill} &
\multicolumn{2}{c@{}}{\hrulefill} \\
$\bolds{k_{0n}}$ & \textbf{1}&
$\bolds{\log p_n}$ & \textbf{1} &$\bolds{\log p_n}$
& \textbf{1} & $\bolds{\log p_n}$ &
\textbf{1} & $\bolds{\log p_n}$ \\
\hline\eqref{eq:prior_loadings1} & 1.73 (1.26) &5.30 (3.92) & 11.92 (2.82) &
13.41 (4.03) &17.01 (0.22) & 19.37 (1.74) & 10.04 (0.08) & 22.10 (0.52) \\
\eqref{eq:shprior_prop} &2.44 (1.40) & 5.42 (2.67) & \phantom{0}4.12 (2.86) &
\phantom{0}7.98 (3.23) & \phantom{0}2.01 (1.44) & \phantom{0}4.56 (1.49) &
\phantom{0}2.04 (1.12) & \phantom{0}5.23 (2.10) \\[3pt]
POET & 3.59 (0.84) & 7.16 (1.84) & \phantom{0}7.14 (1.59) &12.63 (2.89) &
\phantom{0}3.93 (1.17)
& \phantom{0}7.39 (1.69) &\phantom{0}3.90 (0.71) &10.13 (2.09) \\
AT & 2.32 (1.49) & 5.50 (3.09) & \phantom{0}4.04 (2.99) &\phantom{0}8.26 (4.16) &
\phantom{0}2.12 (1.62) &
\phantom{0}4.45 (1.63) & \phantom{0}1.97 (0.90) & \phantom{0}4.96 (2.28) \\
SC & 3.63 (0.88) & 7.32 (1.95) & \phantom{0}7.26 (1.66) &12.85 (3.07) &
\phantom{0}3.95 (1.19) &
\phantom{0}7.44 (1.75)& \phantom{0}3.88 (0.72) & 10.24 (0.26) \\
\hline
\end{tabular*}}
%
%
\end{table}

The results for \eqref{eq:shprior_prop} and \eqref{eq:prior_loadings1}
were reported based on 10,000 runs of the Gibbs sampler with 5000
burn-in. From Tables~\ref{tab:1} and \ref{tab:2}, it becomes evident
that when the number of model parameters increase, the performance of
\eqref{eq:prior_loadings1} deteriorates due to possibly slower
convergence of the MCMC, while (\ref{eq:shprior_prop}) has more robust
performance. Even in Table~\ref{tab:2}, where the truth is misspecified
for both \eqref{eq:shprior_prop} and AT, and in fact designed to favor
POET, \eqref{eq:shprior_prop} performs at least equally or better than
the competitors. For each MCMC iteration, the runtime for \eqref
{eq:shprior_prop} scaled approximately linearly with $n$ and $p$,
though we are not aware of sharp theoretical bounds on MCMC convergence
in high dimensions guaranteeing polynomial time convergence unlike many
frequentist estimators.

\section{Some properties of shrinkage priors in high-dimensional
settings}\label{sec:disc}
We develop a number of properties of the proposed shrinkage prior
\eqref
{eq:shprior_prop} in high-dimensional settings; the results are used to
prove the main results on posterior concentration, but are also of
independent interest. Proofs of all the results are deferred to the \hyperref[app]{Appendix}.

Let $\theta$ be a $p$-dimensional vector and $\theta_0 \in l_0[s; p]$
be an $s$-sparse vector with $s \ll p$. Depending on the problem,
$\theta$ might correspond to a high-dimensional mean vector, a vector
of regression coefficients or a column of the factor loadings, with
$\theta_0$ corresponding to a sparse truth.\footnote{For us, $\theta$
and $\theta_0$ correspond to the vectorized loadings $\Lambda_n$ and
$\Lambda_{0n}$, respectively.} A quantity of fundamental importance in
studying the behavior of the posterior distribution in high-dimensional
problems is the \emph{prior concentration} around an arbitrary sparse
vector $\theta_0$, which is defined as the noncentered small ball probability
%
\begin{equation}
\label{eq:prior_conc_basic} \bbP\bigl(\Vert\theta- \theta_0\Vert_2 <
\varepsilon\bigr),
\end{equation}
for $\varepsilon$ small. It can be shown that if $\theta_j$'s are i.i.d.
standard normal,
\[
\sup_{\theta_0 \in l_0[s; p]} \bbP\bigl(\Vert\theta- \theta_0
\Vert_2 < \varepsilon \bigr) \leq e^{- C p \log(1/\varepsilon)},
\]
which decays exponentially with $p$ for fixed $s$ limiting the ability
of the posterior to concentrate on sparse $\theta_0$.
However, with appropriate point mass mixture priors having a
probability mass at zero and $\Vert\theta_0\Vert$ bounded, the small ball
probability~(\ref{eq:prior_conc_basic}) can be improved to $e^{- C s
\log(1/\varepsilon)}$ \cite{castilloneedles}.

For reasons mentioned in Section~\ref{sec:prior}, there has been a
recent thrust on developing \emph{one-group} alternatives to the \emph
{two-group} mixture priors using continuous shrinkage priors, which can
be often represented as a \emph{global--local} scale mixture \cite
{polson2010shrink} of Gaussians. Despite computational advantages with
this family of shrinkage priors, their concentration properties are
understudied. Our proposed prior \eqref{eq:shprior_prop}, which can be
expressed as a Gaussian scale mixture, favors a large subset of the
$\theta_j$ to be \emph{simultaneously} close to zero while inflicting
minimal shrinkage on the rest, and thus achieve a concentration similar
to point mass mixture priors. In the following Lemma~\ref
{lem:shrinkage_conc}, we present a \emph{nonasymptotic} bound to the
prior concentration for \eqref{eq:shprior_prop}.

\begin{lem}\label{lem:shrinkage_conc}
Suppose $\theta\sim$ \eqref{eq:shprior_prop}. Let $\theta_0 \in l_0[s;
p], 1\leq s \leq p$ and $s/p \leq1/2$. Then, for any $\varepsilon\in
(0,1)$ small enough,
\[
\bbP\bigl(\Vert\theta- \theta_0\Vert_2 < \varepsilon\bigr) \geq
\exp\bigl[-C \max \bigl\{ \Vert\theta_0\Vert_2^2,
s\log(s/\varepsilon), \log p\bigr\} \bigr]
\]
for some constant $C > 0$.
\end{lem}

We also state an auxiliary Lemma~\ref{lem:doubleconc} which is used to
prove Lemma~\ref{lem:shrinkage_conc}; refer to the supplemental
document for a proof.
%
\begin{lem}\label{lem:doubleconc}
Let $\eta\in\mathbb{R}^s$ denote a random vector with independent
components $\eta_j \sim\mathrm{DE}(\psi_j)$. If there exist numbers
$a, b > 0$, such that $\psi_j \in[a, b]$ for all $j = 1, \ldots, s$,
then for any $\delta> 0$ and $\eta_0 \in\mathbb{R}^s$,
\[
\bbP\bigl(\Vert\eta- \eta_0\Vert_2 < \delta\bigr) \geq\exp
\Biggl\{-\frac
{C_1}{a^2}\sum_{j=1}^s
\vert\eta_{0j}\vert^2 - C_2s - s\bigl\vert\log \bigl
\{ \delta /(b\sqrt{s})\bigr\} \bigr\vert \Biggr\}
\]
for constant $C_1, C_2 > 0$.
\end{lem}
%

We next show that the shrinkage prior \eqref{eq:shprior_prop} does not
spread its mass across too many dimensions. A point mass mixture prior
allows a high-dimensional vector to collapse onto fewer dimensions and
the implied dimensionality can be naturally studied through appropriate
tail bounds for the induced prior on $|\operatorname{supp}(\theta)|$, which is
a random variable supported on $\{0,1, \ldots, p\}$. Such bounds on the
prior dimensionality are useful to control the posterior model size
\cite{castilloneedles}. However, continuous shrinkage priors do not
allow exact zeroes in $\theta$ and clearly $\bbP(|\operatorname{supp}(\theta)|
= p) = 1$. We instead use a generalized definition of the support of a
vector as the subset of entries which are larger than a small number
$\delta$ in magnitude. For any $\delta> 0$, we denote the
corresponding subset to be $\operatorname{supp}_{\delta}(\theta)$, so that
$\operatorname{supp}_{\delta}(\theta) = \{ j \dvtx \vert\theta_j\vert >
\delta\}$.\vadjust{\goodbreak}

In the following Lemma~\ref{lem:shrinkage_dim}, we provide a \emph
{nonasymptotic} tail bound for $|\operatorname{supp}_{\delta}(\theta)|$, the
\emph{number} of entries in $\theta$ larger than $\delta$ in magnitude.
%
\begin{lem}\label{lem:shrinkage_dim}
Let $\varepsilon\in(0, 1)$ and $\delta= \varepsilon/p$ with $\varepsilon>
1/p^B$ for some $B > 0$. If $\theta$ is drawn according to the prior
\eqref{eq:shprior_prop} and $s \gtrsim\log p$, then there exists a
constant $A > 0$ such that
\[
\bbP\bigl(\bigl|\operatorname{supp}_{\delta}(\theta)\bigr| > A s\bigr) \leq
e^{- C s}
\]
for some constant $C > 0$. Moreover, the constant $C$ appearing in the
exponent can be made arbitrarily large by choosing $A$ large enough.
\end{lem}
A final important property of \eqref{eq:shprior_prop} is established
through the following deviation result on the $l_1$ norm of $\theta$.
%
\begin{lem}\label{lem:largedev_ps}
For $t \geq1$, $\bbP[\Vert\theta\Vert_1 \geq t ] \leq2 e^{- C
\sqrt{t}}$.
\end{lem}

\section{Construction of test functions}\label{sec:test}
An important step \cite{ghosal2000convergence} in Bayesian asymptotic
theory for establishing posterior contraction rates is to develop a
test function for the true parameter versus the complement of a ball of
radius $\varepsilon> 0$ (in an appropriate norm) around the truth with
type-I and II error rates of the order $\exp(- C n \varepsilon^2)$. Under
the Hellinger or $L_1$ distance between densities, existence of such
tests is guaranteed by the seminal work of \cite{birge1984tests} and
\cite{le1986asymptotic}; the same is true for norms compatible to the
above norms \cite{ghosal2007noniid}. However, when the object of
interest is \emph{not the density itself}, but rather some
high-dimensional parameter indexing the density with a norm of
discrepancy relevant to the space the parameter lives in, the test
arising from Birg{\'e}--Le Cam theory might fail to produce the
desired error rates in the \emph{norm of interest}.

In the context of nonparametric function estimation in general $L_r$
norms, \cite{gine2011rates} advocated using concentration inequalities
based on empirical process techniques as an alternative to the
traditional testing framework. Castillo and Van Der Vaart~\cite
{castilloneedles} used deviation bounds for the likelihood ratio test
in estimating a high dimensional mean in \emph{Euclidean norm}. An
important contribution of the present paper is to utilize recently
developed concentration results for random (self-adjoint) matrices
\cite{tropp2012user,vershynin2010introduction} to devise a test function.

Using a version of the matrix Bernstein inequality (Theorem~6.2 in
\cite
{tropp2012user}), it can be shown that the sample estimator $\hat
{\Sigma
} = n^{-1} \sum_{i=1}^n y_i y_i^{\T}$ has appropriate concentration
around $\bbE\hat{\Sigma} = \Sigma_0$ when the ``effective rank''
$r_e(\Sigma_0):= \tr(\Sigma_0)/\| \Sigma_0 \|_2$ is modest
compared to
$p$ \cite{vershynin2010introduction,bunea2012sample}. However, for
$\Sigma_0 = \Lambda_0 \Lambda_0^{\T} + \sigma_0^2 \mathrm{I}_p
\in
\mathcal{C}_{0n}$, $r_e(\Sigma_0)$ can scale in the order of $p$,
prohibiting us from using $\hat{\Sigma}$ as an estimator to construct
the test. A crucial observation is that even if $\Sigma_0 \in\mathcal
{C}_{0n}$ does not necessarily have a small effective rank, the larger
contribution to the operator norm of $\Sigma_0$ ($\| \Sigma_0 \|_2 =
\Vert\Lambda_0\Vert_2^2 + \sigma_0^2$) comes from the low rank part by
{(A3)}. We exploit this to design a novel projection based test in
Theorem~\ref{thmm:test_opnorm_basic} below, where the types~I and~II
error rates can be expressed in terms of deviation bounds of a $k_{0n}
\times k_{0n}$ sample covariance matrix from its mean. Dependence of
all quantities on $n$ has been made explicit from this point onwards.


\begin{thmm}\label{thmm:test_opnorm_basic}
Recall the sequences $k_{0n}$ and $c_n$ from Assumptions \ref
{ass:truth} and~\ref{ass:truecov1}, respectively. Let $\Sigma_{0n}
\in
\m C_{0n}$ with the corresponding $\Lambda_{0n} \in\mb R^{p_n \times
k_{0n} }$. Let $B_{j,n} = \{ \Sigma_n \in\m C_n \dvtx j \varepsilon_n \leq
\Vert\Sigma_n - \Sigma_{0n} \Vert_2 < (j+1) \varepsilon_n\}$ denote an
annulus of inner radius $j \varepsilon_n$ and outer radius $(j+1)
\varepsilon
_n$ in operator norm around $\Sigma_{0n}$ for some integer $j > 1$ and
sequence $\varepsilon_n > 0$. Assume $\varepsilon_n \geq c_n k_{0n} \sqrt {\frac{\log k_{0n}}{n}}$ and $\varepsilon_n \log j \leq c_n$. Fix
$\Sigma
_{1n} \in B_{j,n}$ and let $E_{j, n} = \{ \Sigma_n \in B_{j, n} \dvtx \Vert\Sigma_n - \Sigma_{1n}\Vert_2 < j \varepsilon_n/2\}$ denote an
operator norm
ball in $B_{j,n}$ around $\Sigma_{1n}$ of radius $j \varepsilon_n/2$.

Based on $n$ i.i.d. samples $y_1, \ldots, y_n$ from $\mathrm{N}_{p_n}(0,
\Sigma_n)$, consider testing the point null vs. composite alternative
hypothesis
%
\begin{equation}
\label{eq:test_op1a_basic} H_0\dvtx \Sigma_n =
\Sigma_{0n}\quad \mbox{versus}\quad H_1\dvtx \Sigma_n
\in E_{j,n}.
\end{equation}
Define $x_i = (1/c_n) \Lambda_{0n}^{\T} y_i$ and $z_i = \Lambda_{0n}
x_i$ for $i = 1, \ldots, n$, so that $x_i \in\mathbb{R}^{k_{0n}}$ and
$z_i \in\mathbb{R}^{p_n}$. Let $\hat{\Sigma}_y = n^{-1} \sum_{i=1}^n
y_i y_i^{\T}$ and define $\hat{\Sigma}_x = (1/c_n^2) \Lambda
_{0n}^{\T}
\hat{\Sigma}_y \Lambda_{0n}, \hat{\Sigma}_z = \Lambda_{0n} \hat
{\Sigma
}_x \Lambda_{0n}^{\T}$. Let $\phi_{j,n}$ denote a test function for
\eqref{eq:test_op1a_basic} defined as
%
\begin{equation}
\label{eq:test_fn_def} \phi_{j,n} = 1_{  \{\Vert\hat{\Sigma}_z - \Sigma_{0n} \Vert_2
\geq j
\varepsilon_n/4  \}}.
\end{equation}
Then, the type-I and type-II error rates of $\phi_{j,n}$ satisfy:
%
\begin{eqnarray}
\bbE_{0} \phi_{j,n} &\leq&\exp \biggl\{-
\frac{C n j^2 \varepsilon
_n^2}{c_n^2k_{0n}^2} \biggr\}\label{eq:errorprob_opnorm_basic1},
\\
 \sup_{\Sigma_n \in E_{j,n} } \bbE_{\Sigma_n} (1 - \phi_{j,n})
&\leq &\exp \biggl\{- \frac{ C n (\log j)^2 \varepsilon_n^2}{c_n^2k_{0n}^2} \biggr\} \label{eq:errorprob_opnorm_basic2}
\end{eqnarray}
for some constant $C> 0$, where $\bbE_{\Sigma_n}$ denotes an
expectation under the distribution of $\mathbf{y}^{(n)}$ under $\mathrm
{N}(0, \Sigma_n)$ and $\bbE_0$ is a shorthand for $\bbE_{\Sigma_{0n}}$.
\end{thmm}
\begin{remark*} If the condition $\varepsilon_n \log j \leq c_n$ is
replaced by $j^{\delta} \varepsilon_n \leq c_n$ for some $0 < \delta
\leq
1$, the type-II error bound in \eqref{eq:errorprob_opnorm_basic2}
becomes $\exp\{- C n j^{2 \delta} \varepsilon_n^2/(c_n^2 k_{0n}^2)\}$.
\end{remark*}
\begin{pf*}{Proof of Theorem \ref{thmm:test_opnorm_basic}}
We shall make use of a matrix concentration result from \cite
{bunea2012sample}. Let $u_1, \ldots, u_n \stackrel{\mathrm{i.i.d.}}
\sim
\mathrm{N}_q(0, \Sigma)$ and $\hat{\Sigma}:= n^{-1} \sum_{i=1}^n u_i
u_i^{\T}$ denote the sample covariance matrix. Proposition A.4 in
\cite
{bunea2012sample} implies that for any $s > 0$ such that $s + \log q < n$,
%
\begin{equation}
\label{eq:colt} \bbP \biggl[ \| \hat{\Sigma} - \Sigma\|_2 > C \tr(
\Sigma) \sqrt {\frac{s
+ \log q}{n}} \biggr] \leq e^{-s}.
\end{equation}

We adapt a fact from Lemma~5.36 of \cite{vershynin2010introduction}.
For a $p \times k$ matrix $B$ with $p > k$, suppose
$\Vert B^{\T} B - \mathrm{I}_k\Vert_2 \leq\max\{\delta, \delta
^2\}$
for some $\delta> 0$. Then
%
\begin{equation}
\label{lem:vershynin} 1 - \delta\leq s_{\min}(B) \leq s_{\max}(B)
\leq1 + \delta.
\end{equation}

Finally, we index matrices that appear frequently in the sequel. Define
%
%
\begin{eqnarray}
\label{eq:psi_n} G_n&: =& \frac{1}{c_n} \Lambda_{0n}^{\T}
\Lambda_{0n} - \mathrm {I}_{k_{0n}},\qquad \Psi_n:=
\frac{1}{c_n} \Lambda_{0n} \Lambda _{0n}^{\T}
- \mathrm{I}_{p_n},
\nonumber
\\[-8pt]
\\[-8pt]
\nonumber
 \Gamma_n&:=& \biggl(1 +
\frac{\sigma
_{0n}^2}{c_n} \biggr) \mathrm{I}_{k_{0n}}.
\end{eqnarray}
Note that $G_n$ is the matrix appearing in {(A3)}. The nonzero
eigenvalues of $G_n$ and $\Psi_n$ are the same; hence, $\Vert G_n\Vert
_2 =
\Vert\Psi_n\Vert_2$.

\textit{Type-I error}: Recall $\hat{\Sigma}_z$ and $\hat
{\Sigma
}_x$ from the theorem statement. We proceed to bound
$\bbE_0 \phi_{j,n} = \bbP_0  [\Vert\hat{\Sigma}_z - \Sigma
_{0n}\Vert_2
\geq j \varepsilon_n/4  ]$.
By the triangle inequality and Lemma~1.1 in the supplemental document,
%
\begin{eqnarray}
\label{eq:dev_op} && \Vert\hat{\Sigma}_z - \Sigma_{0n}
\Vert_2\nonumber \\
&&\qquad\leq\biggl\Vert\Lambda _{0n} \hat {\Sigma}_x
\Lambda_{0n}^{\T} - \Lambda_{0n} \Lambda
_{0n}^{\T} - \frac {\sigma_{0n}^2}{c_n} \Lambda_{0n}
\Lambda _{0n}^{\T} \biggr\Vert_2 +
\sigma_{0n}^2 \Vert\Psi_n\Vert_2
\\
&&\qquad \leq\Vert\Lambda_{0n}\Vert_2^2 \Vert\hat{
\Sigma}_x - \bbE_0 \hat {\Sigma }_x
\Vert_2 + \Vert\Lambda_{0n}\Vert_2^2
\Vert\bbE _0 \hat{\Sigma}_x - \Gamma _n
\Vert_2 + \sigma_{0n}^2 \Vert G_n
\Vert_2.\nonumber
\end{eqnarray}
A simple calculation yields
%
%
\begin{eqnarray}
\label{eq:E0Sx} \bbE_0 \hat{\Sigma}_x & =&
\frac{1}{c_n^2} \Lambda_{0n}^{\T} \bigl[\Lambda
_{0n} \Lambda_{0n}^{\T} + \sigma_{0n}^2
\mathrm{I}_{p_n}\bigr] \Lambda_{0n}
\nonumber
\\[-8pt]
\\[-8pt]
\nonumber
& =& \biggl(
\frac{1}{c_n} \Lambda_{0n}^{\T} \Lambda_{0n}
\biggr)^2 + \frac
{\sigma_{0n}^2}{c_n} \frac{1}{c_n} \Lambda_{0n}^{\T}
\Lambda_{0n}.
\end{eqnarray}
Substituting this in \eqref{eq:dev_op} and using triangle inequality,
the sum of the second and third term in \eqref{eq:dev_op} can be
bounded above by
\[
\Vert\Lambda_{0n}\Vert_2^2 \biggl\Vert \biggl(
\frac{1}{c_n} \Lambda _{0n}^{\T} \Lambda_{0n}
\biggr)^2 - \mathrm{I}_{k_{0n}} \biggr\Vert_2 + \sigma
_{0n}^2(1/c_n + 1) \Vert G_n
\Vert_2.
\]
Recall $\Vert G_n\Vert_2 = o(k_{0n} \sqrt{\log k_{0n}/n})$ by {
(A3)}. In
addition, $\Vert\Lambda_{0n}\Vert_2 \leq2 \sqrt{c_n}$ by {
(A3)} and
$\sigma_{0n}^2 \leq c_n$ from {(A4)}. Note that $\Vert A -
\mathrm {I}_{k_{0n}}\Vert_2 < \delta$ for $A$ symmetric and some
$\delta\in(0,1)$
implies that $\Vert A^2 - \mathrm{I}_{k_{0n}}\Vert_2 \leq3 \delta$. Using
these facts, the expression in the above display can be bounded above
by $(13 c_n + 1) \Vert G_n\Vert_2$. Since we have assumed $\varepsilon_n
\geq
c_n k_{0n} \sqrt{\frac{\log k_{0n}}{n}}$ in the condition of the
theorem, $(13 c_n + 1) \Vert G_n\Vert_2$ can be bounded above by $j
\varepsilon
_n/8$ for $n$ large enough. Substituting this bound in \eqref{eq:dev_op},
\[
\bbP_0 \bigl[\Vert\hat{\Sigma}_z -
\Sigma_{0n}\Vert_2 \geq j \varepsilon_n/4 \bigr]
\leq\bbP_0 \bigl[\Vert\Lambda_{0n}\Vert_2^2
\Vert\hat {\Sigma}_x - \bbE_0 \hat{
\Sigma}_x\Vert_2 \geq j \varepsilon_n/8
\bigr].
\]
Using $\Vert\Lambda_{0n}\Vert_2^2 \lesssim c_n$ one more time, we have
%
\begin{equation}
\bbE_0 \phi_{j,n} \leq\bbP_0 \bigl[ \Vert\hat{
\Sigma}_x - \bbE _0 \hat {\Sigma}_x
\Vert_2 \geq C j \varepsilon_n/c_n \bigr].
\label
{eq:err_rates1}
\end{equation}
By definition, $\hat{\Sigma}_x = n^{-1} \sum_{i=1}^n x_i x_i^{\T}$
is a
$k_{0n} \times k_{0n}$ sample covariance matrix. We now invoke \eqref
{eq:colt} to bound the deviation of $\hat{\Sigma}_x$ from its
expectation under $\bbP_0$
in \eqref{eq:err_rates1}. From \eqref{eq:E0Sx} and using $\tr(A^2) =
\Vert A\Vert_F^2$ for $A$ symmetric, $\tr(\bbE_0 \hat{\Sigma}_x) =
\Vert\Lambda_{0n}^{\T} \Lambda_{0n}/c_n\Vert_F^2 + (\sigma
_{0n}^2/c_n) \Vert\Lambda_{0n}/\sqrt{c_n}\Vert_F^2$. Recall $\sigma
_{0n}^2 \leq c_n$ by
{(A4)}. By Lemma~1.1 in the supplemental document, $\Vert\Lambda
_{0n}^{\T} \Lambda_{0n}/c_n\Vert_F \leq\Vert\Lambda_{0n}/\sqrt {c_n}\Vert_F
\Vert\Lambda_{0n}/\sqrt{c_n}\Vert_2 \leq2 \Vert\Lambda
_{0n}/\sqrt {c_n}\Vert_F$. Hence, $\tr(\bbE_0 \hat{\Sigma}_x)
\leq5 \Vert\Lambda _{0n}/\sqrt{c_n}\Vert_F^2$. Using $\Vert
\Lambda_{0n}\Vert_F \leq\sqrt{k_{0n}}
\Vert\Lambda_{0n}\Vert_2 \leq2 \sqrt{c_n k_{0n}}$, $\tr(\bbE_0
\hat
{\Sigma}_x)$ can be bounded above by~$C k_{0n}$.

Choose $s = C n j^2 \varepsilon_n^2/(k_{0n}^2 c_n^2)$. Since $\varepsilon_n
\geq c_n k_{0n} \sqrt{\log k_{0n}/n}$, we have $s \gtrsim\log k_{0n}$
and hence $C \tr(\bbE_0 \hat{\Sigma}_x) \sqrt{(s + \log k_{0n})/n}
\lesssim C \tr(\bbE_0 \hat{\Sigma}_x) \sqrt{s/n} \leq C j \varepsilon
_n/c_n$. By \eqref{eq:colt}, the expression in the right-hand side of
\eqref{eq:err_rates1} is then bounded above by $e^{-s} = \exp\{-C n j^2
\varepsilon_n^2/(k_{0n}^2 c_n^2)\}$, proving \eqref{eq:errorprob_opnorm_basic1}.

\textit{Type-II error}: Fix $\Sigma_n \in E_{j,n}$. We proceed
to bound $\bbE_{\Sigma_n} (1 - \phi_{j,n}) = \bbP_{\Sigma_n}
[\Vert\hat{\Sigma}_z - \Sigma_{0n}\Vert_2 < j \varepsilon_n/4
]$. By
repeatedly using the triangle inequality, we obtain
\begin{eqnarray*}
 \Vert\hat{\Sigma}_z - \Sigma_{0n} \Vert_2
&\geq&\Vert\Lambda _{0n}\Vert_2^2\biggl \Vert\hat{
\Sigma}_x - \mathrm{I}_{k_{0n}} - \frac{\sigma
_{0n}^2}{c_n}
\mathrm{I}_{k_{0n}}\biggr\Vert_2 - \sigma_{0n}^2
\biggl\Vert\frac
{1}{c_n} \Lambda_{0n} \Lambda_{0n}^{\T}
- \mathrm{I}_{p_n}\biggr\Vert_2
\\
& \geq&\Vert\Lambda_{0n}\Vert_2^2 \bigl\{\Vert
\bbE_{\Sigma_n} \hat {\Sigma }_x - \Gamma_n
\Vert_2 - \Vert\hat{\Sigma}_x - \bbE _{\Sigma_n}
\hat {\Sigma }_x\Vert_2 \bigr\} - \sigma_{0n}^2
\Vert \Psi_n\Vert_2.
\end{eqnarray*}
Recall $\Vert\Psi_n\Vert_2 = \Vert G_n\Vert_2$. Therefore, on the set
$
\{ \Vert\hat{\Sigma}_z - \Sigma_{0n}\Vert_2 < j \varepsilon_n/4
\}$,
%
\begin{eqnarray} \label{eq:dev_op2}
 \Vert\Lambda_{0n}\Vert_2^2 \Vert\hat{
\Sigma}_x - \bbE_{\Sigma
_n} \hat {\Sigma}_x
\Vert_2& \geq&\Vert\Lambda_{0n}\Vert_2^2
\Vert \bbE_{\Sigma_n} \hat {\Sigma}_x - \Gamma_n
\Vert_2 - \sigma _{0n}^2 \Vert G_n
\Vert_2 - \frac{j
\varepsilon_n}{4}
\nonumber
\\
& \geq&\Vert\Lambda_{0n}\Vert_2^2 \Vert
\bbE_{\Sigma_n} \hat {\Sigma }_x - \bbE_0 \hat{
\Sigma}_x\Vert_2 \\
&&{}- \Vert\Lambda _{0n}
\Vert_2^2 \Vert\bbE _0\hat {
\Sigma}_x - \Gamma_n \Vert_2 -
\sigma_{0n}^2 \Vert G_n\Vert_2 -
\frac
{j\varepsilon
_n}{4}.\nonumber
\end{eqnarray}
Recalling the definition of $\hat{\Sigma}_x$ and invoking Lemma~1.1 in
the supplemental document,
%
\begin{eqnarray}
\label{eq:expdiff}  \Vert\Lambda_{0n}\Vert_2^2
\Vert\bbE_{\Sigma_n} \hat{\Sigma }_x - \bbE_0
\hat{\Sigma}_x\Vert_2 
&=& \biggl\Vert
\frac{\Lambda_{0n}}{\sqrt{c_n}}\biggr\Vert_2^2 \biggl\Vert\frac
{1}{c_n}
\Lambda _{0n}^{\T} (\Sigma_n -
\Sigma_{0n}) \Lambda _{0n}\biggr\Vert_2
\nonumber
\\
& \geq&\biggl\Vert\frac{\Lambda_{0n}}{\sqrt{c_n}}\biggr\Vert_2^2 s_{\min
}
\biggl(\frac
{\Lambda_{0n}^{\T} \Lambda_{0n}}{c_n} \biggr) \Vert\Sigma_n - \Sigma
_{0n}\Vert_2 \\
&\geq&\biggl\Vert\frac{\Lambda_{0n}}{\sqrt {c_n}}
\biggr\Vert_2^2 s_{\min
} \biggl(\frac{\Lambda_{0n}^{\T} \Lambda_{0n}}{c_n}
\biggr) \frac{j
\varepsilon_n}{2}.\nonumber
\end{eqnarray}
The last inequality in \eqref{eq:expdiff} used the triangle inequality
to obtain $\Vert\Sigma_n - \Sigma_{0n}\Vert_2 \geq\Vert\Sigma
_{1n} - \Sigma_{0n}\Vert_2 - \Vert\Sigma_n - \Sigma_{1n}\Vert_2
\geq j \varepsilon_n
- j
\varepsilon_n/2 = j \varepsilon_n/2$.
By {(A3)} and \eqref{lem:vershynin}, both $\Vert\Lambda
_{0n}/\sqrt {c_n}\Vert_2$ and $s_{\min} (\Lambda_{0n}^{\T}
\Lambda_{0n}/c_n )$
can be bounded below by $4/5$. Hence, $\Vert\Lambda_{0n}\Vert_2^2
\Vert\hat{\Sigma}_x - \bbE_{\Sigma_n} \hat{\Sigma}_x\Vert_2$ is bounded
below by
$64 j \varepsilon_n/250$. Further, based on the calculations following
\eqref{eq:dev_op},
\[
\Vert\Lambda_{0n}\Vert_2^2 \Vert
\bbE_0\hat{\Sigma}_x - \Gamma _n
\Vert_2 + \sigma_{0n}^2 \Vert G_n
\Vert_2 +\frac{j \varepsilon_n}{4}
\]
can be bounded above $63 j\varepsilon_n/250$ for $n$ large enough.
Substituting in \eqref{eq:dev_op2},
\[
\bbP_{\Sigma_n}\bigl [\Vert\hat{\Sigma}_z -
\Sigma_{0n}\Vert_2 \leq j \varepsilon_n/2 \bigr] <
\bbP_{\Sigma_n} \bigl[\Vert\Lambda _{0n}\Vert_2^2
\Vert\hat{\Sigma}_x - \bbE_{\Sigma_n} \hat{
\Sigma}_x\Vert_2 > C j \varepsilon _n
\bigr].
\]
As in case of the type-I error, using $\Vert\Lambda_{0n}\Vert_2^2
\lesssim
c_n$, we conclude that
%
\begin{equation}
\label{eq:err_rates2} \bbE_{\Sigma_n}(1 - \phi_{j,n}) \leq
\bbP_{\Sigma_n} \bigl[ \Vert \hat {\Sigma}_x - \bbE_{\Sigma_n}
\hat{\Sigma}_x\Vert_2 > C j\varepsilon
_n/c_n \bigr].
\end{equation}
We are now in a position to invoke \eqref{eq:colt} to bound the
right-hand side of \eqref{eq:err_rates2}.
Using triangle inequality and von Neumann's trace inequality \cite
{mirsky1975trace},\hskip.2pt\footnote{$|\tr(A)| \leq k \|A\|_2$ for a $k \times
k$ matrix $A$.}
\begin{eqnarray*}
\tr(\bbE_{\Sigma_n} \hat{\Sigma}_x)& \leq&\tr(
\bbE_0 \hat{\Sigma }_x) +\bigl | \tr(\bbE_{\Sigma_n}
\hat{\Sigma}_x - \bbE_0 \hat{\Sigma}_x)
\bigr| \\
&\leq& \tr(\bbE_0 \hat{\Sigma}_x) + k_{0n}
\Vert\bbE_{\Sigma_n} \hat {\Sigma }_x - \bbE_0
\hat{\Sigma}_x\Vert_2.
\end{eqnarray*}
Since $\Sigma_n \in B_{j,n}$, $\Vert\Sigma_n - \Sigma_{0n}\Vert_2
\leq
(j+1) \varepsilon_n < 2 j \varepsilon_n$, and hence
\begin{eqnarray*}
\Vert\bbE_{\Sigma_n} \hat{\Sigma}_x - \bbE_0
\hat{\Sigma }_x\Vert_2& =& \bigl\Vert\Lambda_{0n}^{\T}
(\Sigma_n - \Sigma_{0n}) \Lambda _{0n}
\bigr\Vert_2/c_n^2 \\
&\leq &C \Vert\Sigma_n
- \Sigma_{0n}\Vert_2/c_n \leq2C j
\varepsilon_n/c_n.
\end{eqnarray*}
Substituting in the previous display and using $\tr(\bbE_0 \hat
{\Sigma
}_x) \leq C k_{0n}$, one has $\tr(\bbE_{\Sigma_n} \hat{\Sigma}_x)
\lesssim k_{0n} \max\{1, j \varepsilon_n/c_n \} \leq C k_{0n} (j/\log j)$,
with the last inequality using $\varepsilon_n \log j \leq c_n$.

Choosing $s = C n (\log j)^2 \varepsilon_n^2/(k_{0n}^2 c_n^2)$, one has $s
\gtrsim\log k_{0n}$ and $C \tr(\bbE_{\Sigma_n}\hat{\Sigma}_x)\times
\sqrt{(s + \log k_{0n})/n} \leq C j \varepsilon_n/c_n$. By \eqref
{eq:colt}, the expression in the right-hand side of~\eqref
{eq:err_rates2} is bounded above by $e^{-s} = \exp\{ - C n (\log j)^2
\varepsilon_n^2/(k_{0n}^2 c_n^2)\}$. Since the bound is independent of
$\Sigma_n \in E_{j,n}$, \eqref{eq:errorprob_opnorm_basic2} follows.
\end{pf*}

\section{Proof of the main results}\label{sec:main_pf}

We now proceed to prove the results stated in Section~\ref
{sec:main_res}. We prove Theorem~\ref{thmm:k_0} with the shrinkage prior
\eqref{eq:shprior_prop}; the special case of $k_{0n} = O(1)$ in Theorem~\ref{thmm:op_norm_rates_shrinkage} follows immediately. For the point
mass prior, we only sketch an argument. We introduce a number of
auxiliary Lemmata \ref{lem:denom}, \ref{lem:cov_priorconc}, \ref
{lem:support1} whose proofs can be found in the supplemental document.
\subsection{Proof of Theorem \texorpdfstring{\protect\ref{thmm:k_0}}{5.3}}
Set $\varepsilon_n = c_n k_{0n}^{3/2} \sqrt{\frac{s_n \log p_n}{n}}
\sqrt{\log n}$ and define
$U_{n} = \{\Sigma_n\dvtx \Vert\Sigma_{n} -\Sigma_{0n}\Vert_{2} \leq M
\varepsilon
_n\}$.
The posterior probability assigned to the complement of $U_n$ is given by
%
\begin{equation}
\label{eqn:fraceq}  \Pi_n\bigl(U_n^c \mid
\mathbf{y}^{(n)} \bigr) = \frac{\int_{U_n^c} \prod_{i=1}^{n}
 ({f_{\Sigma_n}(y_i)}/{f_{\Sigma_{0n}(y_i)}}) \,d\Pi
_n(\Sigma
_n)}{\int\prod_{i=1}^{n} ({f_{\Sigma_n}(y_i)}/{f_{\Sigma
_{0n}(y_i)}}) \,d\Pi_n(\Sigma_n)} \equiv\frac{\mathcal{N}_n}{\mathcal{D}_n},
\end{equation}
where $f_{\Sigma_n}$ denotes a $p_n$-dimensional $\mathrm{N}(0,\Sigma_n)$
distribution and $\mathcal{N}_n$ and $\mathcal{D}_n$ denote the
numerator and denominator of the
fraction in \eqref{eqn:fraceq}.

Let $\sigma(y_1, \ldots, y_n)$ denote the $\sigma$-field generated by
$y_1, \ldots, y_n$. We first claim that we can lower-bound $\mathcal
D_n$ on an event $A_n \in\sigma(y_1, \ldots, y_n)$ with large
probability under $f_{\Sigma_{0n}}$ in Lemma~\ref{lem:denom}.

\begin{lem}\label{lem:denom}
Let $\Sigma_{0n} \in\mathcal{C}_{0n}$. Let $\eta_n$ be a sequence
satisfying $\eta_n/s_{\min}(\Sigma_{0n}) \to0$ and $n \eta_n^2
/s_{\min
}(\Sigma_{0n})^2 \to\infty$, and define $\varrho_n = 2s_{\max
}(\Sigma
_{0n})/s_{\min}(\Sigma_{0n})$. Then there exists $A_n \in\sigma(y_1,
y_2, \ldots, y_n)$ with $\bbP_{\Sigma_{0n}}(A_n) \to1$ such that on $A_n$,
\[
\mathcal D_n 
\geq e^{-C n \eta_n^2 \log(\varrho_n)/s_{\min}(\Sigma_{0n})^2}
\Pi_n\bigl( \Sigma_n\dvtx \Vert\Sigma_n -
\Sigma_{0n}\Vert_{F} < \eta_n\bigr).
\]
\end{lem}
We shall set $\eta_n = \sqrt{s_n k_{0n}/n}$ in all future usage of
Lemma~\ref{lem:denom}. Based on our prior specification, $\Sigma_n$ can
be parameterized by $(k, \Lambda_n, \sigma_n^2)$ with $k \in\{1,
\ldots, \infty\}, \Lambda_n \in\mb R^{p_n \times k}, \sigma_n^2 \in(0,
\infty)$, and $\Sigma_n = \Lambda_n \Lambda_n^{\T} + \sigma_n^2
\mathrm
{I}_{p_n}$. We use this to bound $\Pi_n(\Vert\Sigma_n - \Sigma
_{0n}\Vert_F
\leq\eta_n)$ from below in the following Lemma~\ref{lem:cov_priorconc}.
%
\begin{lem}\label{lem:cov_priorconc}
If $\Sigma_{0n} \in\m C_{0n}$, the prior $\Pi_n$ on $\Sigma_n$ is as
in Theorem~\ref{thmm:op_norm_rates_shrinkage}, and $\eta_n = \sqrt{s_n
k_{0n}/n}$, then
\[
\Pi_n\bigl(\Vert\Sigma_n - \Sigma_{0n}
\Vert_F \leq\eta_n \bigr) \geq e^{-
C s_n
k_{0n} \log n }.
\]
\end{lem}
We now introduce some notation. Let
%
%
\begin{eqnarray}
\label{eq:deltan} e_n &=& s_n k_{0n} \log n,\qquad
t_n = C e_n^2,\qquad \delta_n =
\varepsilon _n/(e_n t_n),
\nonumber
\\[-8pt]
\\[-8pt]
\nonumber
\delta_n' &= &\delta_n/(p_n
e_n).
\end{eqnarray}
Recalling that the true loadings has $s_n k_{0n}$ many nonzero entries,
$e_n$ can be thought of as an \emph{effective sparsity} parameter. Also,
recall the $\operatorname{supp}_{\delta}$ notation from Section~\ref{sec:disc}.
Given $k$, let $\operatorname{supp}_{\delta_n'}(\Lambda_n)$ denote the set $S
\subset\{1, \ldots, p_n k\}$ corresponding to the entries in $\operatorname
{vec}(\Lambda_n)$ larger than $\delta_n'$ in absolute magnitude.

Since $\bbP_{\Sigma_{0n}}(A_n) \to1$ by Lemma~\ref{lem:denom}, it is
enough to show
\[
\lim_{n \rightarrow\infty} \bbE_{0} \bigl[\Pi_n
\bigl(U_n^c \mid\mathbf {y}^{(n)}
\bigr)1_{A_n} \bigr] = 0
\]
to prove Theorem~\ref{thmm:op_norm_rates_shrinkage}, where $\bbE_0$ is a
shorthand for $\bbE_{\Sigma_{0n}}$.
For some $H > 0$ to be chosen later,
%
\begin{eqnarray}
\label{post:o} \bbE_{0} \bigl[\Pi_n\bigl(U_n^c
\mid\mathbf{y}^{(n)}\bigr)1_{A_n} \bigr] &\leq& \bbE
_{0} \bigl[\Pi_n\bigl(U_n^* \mid
\mathbf{y}^{(n)}\bigr)1_{A_n} \bigr]
\nonumber
\\
&&{}+ \bbE_{0} \bigl[ \Pi_n \bigl(W_n^c
\cap V_n \mid\mathbf {y}^{(n)} \bigr) 1_{A_n}
\bigr] \\
&&{}+ 2 \bbE_{0} \bigl[ \Pi_n\bigl(V_n^c
\mid\mathbf{y}^{(n)}\bigr) 1_{A_n} \bigr],\nonumber \phantom{}
\end{eqnarray}
where
$U_n^* = U_n^c \cap W_n \cap V_n$, with
%
%
\begin{eqnarray}
\label{eq:vnwn} W_n &=& \bigl\{\bigl|\operatorname{supp}_{\delta_n'}(
\Lambda_n)\bigr| \leq H e_n, \Vert \Lambda_n
\Vert_1 \leq t_n, \sigma^2 \leq
t_n \bigr\},
\nonumber
\\[-8pt]
\\[-8pt]
\nonumber
 V_n& =& \{k \leq C e_n \}.
\end{eqnarray}
Thus, $U_n^*$ consists of (strictly speaking, can be identified with
the class of) covariance matrices $\Sigma_n = \Lambda_n \Lambda
_n^{\T}
+ \sigma_n^2 \mathrm{I}_{p_n}$ satisfying $\Vert\Sigma_n - \Sigma
_{0n}\Vert_2 > M \varepsilon_n$, where $\Lambda_n \in\mb R^{p_n \times k}$
with $k \leq C e_n$, $ \Vert\Lambda_n\Vert_1 \leq t_n$, $|\operatorname
{supp}_{\delta_n'}(\Lambda_n)| \leq H e_n$ and $\sigma_n^2 \leq t_n$.

We now show in Lemma~\ref{lem:support1} that the expression in \eqref
{post:o} goes to zero, so that we can focus on $\Pi_n(U_n^* \mid
\mathbf
{y}^{(n)})$. This will be crucial in reducing the entropy of the model space.
%
\begin{lem}\label{lem:support1}
Recall the sequences and sets in \eqref{eq:deltan} and \eqref{eq:vnwn},
respectively. There exist constants $H, C > 0$ such that
\begin{eqnarray*}
 \lim_{n \rightarrow\infty} \bbE_{0} \bigl[ \Pi_n
\bigl(W_n^c \cap V_n \mid
\mathbf{y}^{(n)}\bigr) 1_{A_n} \bigr]& =& 0,
\\
 \lim_{n \rightarrow\infty} \bbE_{0} \bigl[ \Pi_n
\bigl( V_n^c \mid \mathbf {y}^{(n)}\bigr)
1_{A_n} \bigr] &=& 0.
\end{eqnarray*}
\end{lem}


For $k \leq C e_n$, a set $S \subset\{1, \ldots, p_n k\}$ with $|S|
\leq H e_n$, and $j \geq M$, let $B_{k,S,j,n}$ denote the following
subset of $U_n^*$:
%
%
\begin{eqnarray}\qquad
\label{eq:bjsn} B_{k,S,j,n} &=& \bigl\{ \Sigma_n =
\Lambda_n \Lambda_n^{\T} + \sigma_n^2
\mathrm{I}_{p_n} \dvtx \Lambda_n \in\mb R^{p_n \times k}, k
\leq C e_n, \Vert\Lambda_n\Vert_1 \leq
t_n,
\nonumber
\\[-8pt]
\\[-8pt]
\nonumber
&&\hspace*{8pt}{}\sigma_n^2 \leq t_n, \operatorname {supp}_{\delta_n'}(\Lambda_n) = S, j
\varepsilon_n \leq\Vert\Sigma _{n} - \Sigma_{0n}
\Vert_{2} < (j+1) \varepsilon_n \bigr\}.
\end{eqnarray}
Then, using a standard testing argument (see, e.g., the proof of
Proposition~5.1 in~\cite{castilloneedles}),
%
\begin{eqnarray}
\label{eq:main_post} 
&&\bbE_{0} \bigl[
\bbP\bigl(U_n^* \mid\mathbf{y}^{(n)} \bigr) 1_{A_n}
\bigr]
\nonumber
\\
 &&\qquad\leq\sum_{k \leq C e_n }\sum
_{S \dvtx |S| \leq H e_n } \sum_{j \geq M} \Bigl[
\bbE_{0} \Phi_{k,S,j,n}\\
 &&\hspace*{125pt}{}+ \beta_{k,S,j,n} \sup
_{\Sigma_n
\in
B_{k,S,j,n}} \bbE_{\Sigma_{n}} (1- \Phi_{k,S,j,n}) \Bigr],\nonumber
\end{eqnarray}
where $\Phi_{k,S,j,n}$ is a (point vs. composite) test function for
%
\begin{equation}
\label{eq:test_op1} H_0\dvtx \Sigma_n =
\Sigma_{0n} \quad\mbox{versus}\quad H_1\dvtx \Sigma_n
\in B_{k,S,j,n}
\end{equation}
whose construction is provided below
and
%
\begin{equation}
\label{eq:beta_jsn} \beta_{k,S,j,n}:= \frac{\Pi_n(B_{k,S,j,n})}{ e^{-n \eta_n^2 \log
\varrho_n/s_{\min}(\Sigma_{0n})^2} \Pi_n(\Vert\Sigma_n - \Sigma
_{0n}\Vert_F < \eta_n) } \leq e^{C e_n}.
\end{equation}
To obtain the upper bound on $\beta_{k,S,j,n}$ in the above display,
bound $\Pi_n(B_{k,S,j,n})$ above by $1$, use the fact that $\log
\varrho
_n \lesssim\log c_n \leq\log n$ by {(A1)} and {(A4)} and use
Lemma~\ref{lem:cov_priorconc} to conclude that
$\beta_{k,S,j,n} \leq e^{C e_n}$.

To construct the test function $\Phi_{k,S,j,n}$ in \eqref
{eq:main_post}, we cover $B_{k,S,j,n}$ with a union of balls and obtain
local tests for $\Sigma_{0n}$ versus the centers of each of the balls
using Theorem~\ref{thmm:test_opnorm_basic}. Since we are inside $W_n
\cap V_n$, the number of such balls can be controlled and $\Phi
_{k,S,j,n}$ is obtained as the maximum of the local tests.

Let $\Sigma_{n,l}$\footnote{We suppress the dependence on $k, S$ and
$j$.} for $l \in I_{k,S,j,n}$ be a $j \varepsilon_n/2$-net of
$B_{k,S,j,n}$ in \emph{operator norm} and for each $l$, define
$E_{k,S,j,n,l} = \{\Sigma_n \in B_{k,S,j,n} \dvtx \Vert\Sigma_n -
\Sigma _{n,l}\Vert_2 \leq j \varepsilon_n/2\}$. By definition,
\[
B_{k,S,j,n} \subset\bigcup_{l \in I_{k,S,j,n}}
E_{k,S,j,n,l}.
\]
Clearly, $j \varepsilon_n \leq\Vert\Sigma_{n,l} - \Sigma_{0n}\Vert
_{2} <
(j+1) \varepsilon_n$ and $\varepsilon_n \geq c_n k_{0n} \sqrt{\log
k_{0n}/n}$. For $\Sigma_n \in U_n^*$, $\Vert\Lambda_n\Vert_2 \leq
\sqrt {e_n} t_n$ and $\Vert\Sigma_n - \Sigma_{0n}\Vert_2 \leq\Vert
\Sigma_n\Vert_2
+ c_n \leq e_n t_n^2 + t_n + c_n \lesssim e_n^5$. Hence, $B_{k,S,j,n}
\subset U_n^*$ implies $j \lesssim e_n^5/\varepsilon_n$, and hence $\log j
\lesssim\log n$, so that $\varepsilon_n \log j \leq c_n$ by the second
part of {(A1)}.
Therefore, the conditions of Theorem~\ref{thmm:test_opnorm_basic} for
the \emph{point versus composite} test $H_0\dvtx \Sigma_n = \Sigma_{0n}$
versus $H_1\dvtx \Sigma_n \in E_{k,S,j,n,l}$ are satisfied. Let $\phi
_{k,S,j,n,l}$ denote the corresponding test function from Theorem~\ref
{thmm:test_opnorm_basic} with \mbox{type-I} error $e^{-C n j^2 \varepsilon
_n^2/(c_n^2 k_{0n}^2)} = e^{-C j^2 e_n \log p_n}$ and type-II error
$e^{-C n (\log j)^2 \varepsilon_n^2/(c_n^2 k_{0n}^2)} = e^{-C (\log j)^2
e_n \log p_n}$.
Letting $\Phi_{k,S,j,n} =\break  \max_{l \in I_{k,S,j,n}}
\phi_{k,S,j,n,l}$, we therefore have
\begin{eqnarray*}
 \bbE_{0} (\Phi_{k,S,j,n}) &\leq&| I_{k,S,j,n} |
e^{-C j^2 e_n \log
p_n },
\\
\sup_{\Sigma_n \in B_{j,S, n}} \bbE_{\Sigma_n}(1- \Phi_{k,S,j,n})
&\leq& e^{-C (\log j)^2 e_n \log p_n}.
\end{eqnarray*}
%

To estimate $|I_{k,S,j,n}|$, {that is}, the covering number of
$B_{k,S,j,n}$ in operator norm, we embed $B_{k,S,j,n}$ inside $\tilde
{B}_{k,S,j,n}$, whose covering number is easier to calculate:
\begin{eqnarray*}
\tilde{B}_{k,S,j,n}&:=& \bigl\{\Sigma_n =
\Lambda_n \Lambda_n^{\T} + \sigma
_n^2 \mathrm{I}_{p_n}\dvtx
\Lambda_n \in B_{k,S,j,n}^{(\Lambda)}, \sigma_n^2
\leq t_n,
\\
&&\hspace*{61pt}{} j \varepsilon_n \leq\Vert\Sigma_n -
\Sigma_{0n}\Vert_2 < (j+1) \varepsilon _n\bigr
\},
\end{eqnarray*}
where $B_{k,S,j,n}^{(\Lambda)} = \{ \Lambda_n \in\mb R^{p_n \times k}
\dvtx k \leq C e_n, \mathrm{supp}_{\delta_n'}(\Lambda_n) = S, \Vert
\Lambda_n\Vert_F \leq e_n t_n \}$.
The containment $B_{k,S,j,n} \subset\tilde{B}_{k,S,j,n}$ follows since
$\Vert\Lambda_n\Vert_F \leq\sqrt{k} \Vert\Lambda_n\Vert_2 \leq
k \Vert\Lambda_n\Vert_1 \leq e_n t_n$.

We now proceed to explicitly construct a $j \varepsilon_n/2$-net for
$\tilde{B}_{k,S,j,n}$. Let $\xi_n = j \varepsilon_n/(8 e_n t_n)$. For
notational convenience, we use $P_S(\theta)$ below to denote $\theta_S$
defined in Section~\ref{sec:prelim}.
Let $\{\Lambda_l\}_{l=1}^{L}$ be a $\xi_n$-net of
$B_{k,S,j,n}^{(\Lambda
)}$. Also, let $\{ \sigma_r^2 \}_{r=1}^{R}$ be a $j \varepsilon_n/4 $-net
of $[0, t_n]$. We show below that $\{ \Lambda_l \Lambda_l^{\T} +
\sigma
_r^2\}_{l,r}$ form a $j \varepsilon_n/2$-net of $\tilde{B}_{k,S,j,n}$ in
operator norm.

Let $\tilde{\Sigma} = \tilde{\Lambda} \tilde{\Lambda}^{\T} +
\tilde
{\sigma}^2 \mathrm{I}$ be in $\tilde{B}_{k,S,j,n}$. Find $\Lambda_l$
and $\sigma_r^2$ from the respective nets so that $\Vert\Lambda_l -
\tilde{\Lambda}\Vert_F \leq\xi_n$ and $| \sigma_r^2 - \tilde
{\sigma}
^2 |
\leq j \varepsilon_n/4$. Let $\Sigma= \Lambda_l \Lambda_l^{\T} +
\sigma
_r^2$. Then
\[
\Vert\Sigma- \tilde{\Sigma}\Vert_2 \leq j\varepsilon_n/4
+ \bigl\Vert \Lambda_l \Lambda_l^{\T} - \tilde{
\Lambda} \tilde{\Lambda}^{\T
}\bigr\Vert \leq j\varepsilon_n/4 + \bigl[
\Vert \Lambda_l \Vert _2 + \Vert \tilde{\Lambda} \Vert
_2 \bigr] \xi_n \leq j \varepsilon_n/2.
\]
We have thus proved our claim, and hence $|I_{k,S,j,n}| \leq L \times
R$. Note the use of the control on $\Vert\Lambda\Vert_2$ over
$B_{k,S,j,n}$ in the above display.

Clearly, $R$ can be chosen to smaller than $t_n/(2j \varepsilon_n)$. With
$s = |S|$, let $\{ \theta_l \}_{l=1}^{L}$ be a $\xi_n/2$-net of the
Euclidean sphere in $\mathbb{R}^s$ of radius $e_n t_n$. By Lemma~5.2 of
\cite{vershynin2010introduction}, the cardinality of such a net $L$ can
be chosen to be smaller than $(1 + e_n t_n/\xi_n)^s$. We now exhibit a
$\xi_n$-net $\{\Lambda_l\}_{l=1}^{L}$ for $B_{k,S,j,n}^{(\Lambda)}$ in
Frobenius norm as follows. Set $P_S(\Lambda_l) = \theta_l$ and
$P_{S^c}(\Lambda_l) = \mathbf{0}$. Let $\Lambda\in B_{k,S,j,n}^{(\Lambda
)}$ and $\theta= P_S(\Lambda)$.\vadjust{\goodbreak} There exists $\theta_l$ such that
$\Vert\theta- \theta_l\Vert_2 \leq\xi_n/2$. Also, since $\operatorname
{supp}_{\delta_n'}(\Lambda) = S$, $\Vert P_{S^c}(\Lambda) \Vert_2
\leq
\delta_n$. By choosing $j$ larger than some constant $J$, we can make
$\xi_n \geq2 \delta_n$. Hence, $\Vert\Lambda_l - \Lambda\Vert_F
\leq
\xi
_n$. Thus, $L \times R$ can be bounded above by $e^{C s \log(e_n t_n)}
\leq e^{C s \log p_n}$, and hence
%
\begin{eqnarray}
 \bbE_{0} (\Phi_{k,S,j,n})& \leq& e^{C s \log p_n}
e^{-C_1 j^2 e_n
\log
p_n}, \label{eq:error_opnorm1}
\\
 \sup_{\Sigma_n \in B_{k,S,j,n}} \bbE_{\Sigma_n}(1- \Phi_{k,S,j,n})
&\leq& e^{-C_2 (\log j)^2 e_n \log p_n}. \label{eq:error_opnorm2}
\end{eqnarray}

Substitute the bounds obtained in \eqref{eq:beta_jsn}, \eqref
{eq:error_opnorm1} and \eqref{eq:error_opnorm2} in \eqref
{eq:main_post}. Observing that all the bounds are free of $k$, we can
bound the expression in (\ref{eq:main_post}) by
%
\begin{equation}
\label{eq:finalex} (C e_n) \sum_{s=0}^{H e_n}
\pmatrix{p_n \cr s} \biggl[\sum_{j \geq M}
e^{C s
\log p_n} e^{-C_1j^2 e_n \log p_n} + e^{C e_n} e^{-C_2 (\log j)^2 e_n
\log p_n} \biggr].\hspace*{-45pt}
\end{equation}
The first term in the inner sum over $j$ can be bounded above by $e^{C
e_n \log p_n} \times e^{- C_3 M^2 e_n \log p_n}$, while the second one
by $e^{C e_n} e^{- C_4 (\log M)^2 e_n \log p_n}$. Noting that $(C e_n)
\times(H e_n + 1) \max_{\{0\leq l \leq H e_n\}} {p_n \choose l} \leq
\exp\{C e_n \log p_n\}$, (\ref{eq:finalex}) goes to $0$ as $n \to
\infty
$ for a large enough constant $M > 0$. This completes the proof of
Theorem~\ref{thmm:k_0} with the shrinkage prior \eqref{eq:shprior_prop}.

The proof for the point mass priors \eqref{eq:prior_loadings1} follows
similarly. Since the point mass mixture priors allow exact zeros in the
loadings, we can condition on $\operatorname{supp}(\Lambda_n) = S$. By
properties of point mass mixture priors shown in \cite
{castilloneedles}, analogues of Lemmata~\ref{lem:cov_priorconc} and
\ref
{lem:support1} can be obtained to complete the theorem.


\begin{appendix}\label{app}
\section*{Appendix}
\begin{pf*}{Proof of Lemma \protect\ref{lem:fan1}}
Observe that if $\| \theta_j - \theta_{j'} \|_H = 2(s-r)$, then
$\langle b_j, b_{j'}\rangle = r$. For $j \neq j'$, $\Sigma_{(j)} -
\Sigma_{(j')} =
\gamma(b_jb_j^{\T} - b_{j'} b_{j'}^{\T})$.
The nonzero eigenvalues of the matrix $B = (b_j b_j^{\T} - b_{j'}
b_{j'}^{\T})$
are $\{\sqrt{s^2 - r^2},-\sqrt{s^2 - r^2} \}$, since $\operatorname{rk}(B) =
2$, $\tr(B) = 0$ and $\tr(B^2) = 2(s^2 - r^2)$.

Since $\theta_j \in\m M$ for all $j$, by symmetry, $\mrr{det}(\Sigma
_{(j)}) = \mrr{det}(\Sigma_{(j')})$ for all $j \neq j'$. Hence, $\mrr
{KL}(\bbP^{(j)}, \bbP^{(j')}) = (n/2) \{ \tr( \Sigma_{(j)}^{-1}
\Sigma
_{(j')}) - p \}$. Write $\Sigma_{(j)} = \beta(A + t b_j b_j^{\T})$,
where $A$ is a diagonal matrix with the first $(p-1)$ diagonal entries
equaling one and the $p${th} entry being $(1 + \kappa/\beta)$. An
application of the Woodbury matrix inversion formula produces
\[
\bigl(A + t b_j b_j^{\T}
\bigr)^{-1} = A^{-1} - \frac{t }{1 + ts} b_j
b_j^{\T},
\]
so that
\[
\Sigma_{(j)}^{-1} \Sigma_{(j')} =
\mr{I}_p - \frac{t }{1 + t s} b_j b_j^{\T}
+ t b_{j'} b_{j'}^{\T} - \frac{t^2 r}{1+ t s}
b_j b_{j'}^{\T}.
\]
The proof is completed by observing that $\tr(b_j b_j^{\T}) = s$ and
$\tr(b_j b_{j'}^{\T}) = r$.\vadjust{\goodbreak}
\end{pf*}

\begin{pf*}{Proof of Lemma \protect\ref{lem:pack}} Let $\tau\in\m
M$ with
$\mrr{supp}(\tau) = S$. We show that for any $x \in\m M$, $\| x -
\tau
\|_H = 2 \sum_{j \notin S} 1(x_j = 1)$. To that end, we have $\| x -
\tau\|_H = \sum_{j \in S} 1(x_j = 0) + \sum_{j \notin S} 1(x_j = 1) =
s + a - b$, where $a = \sum_{j \notin S} 1(x_j = 1)$, $b = \sum_{j
\in
S} 1(x_j = 1)$. Since $x \in\m M$, we also have $a + b = s$, which
implies $\|x - \tau\|_H = 2a$.

Let $k$ denote the integer part of $s/6$. Let $\m M_0$ be a maximal set
of points in $\m M$, with each pair at least $2(k+1)$ apart in Hamming
distance. Note here that $2(k+1) > s/3$. Since $\m M_0$ is maximal and
$d(x, y)$ is even for any $x, y \in\m M$ by the above calculation, it
follows that $\m M \subset\bigcup_{\tau\in\m M_0} B(\tau; 2k)$, where
\[
B(\tau; 2k) = \bigl\{ x \in\m M \dvtx \| x - \tau\|_H \leq2k \bigr\}.
\]
By symmetry, $B(\tau; 2k)$ is independent of $\tau$, so that $| \m M |
\leq| \m M_0 | | B(\tau; 2k)|$ for any $\tau\in\m M_0$. It is easy
to see that
\[
\bigl| B(\tau; 2k) \bigr| = \sum_{j=0}^k |
A_j| = \sum_{j=0}^k \pmatrix{s
\cr j} \pmatrix{q-s \cr j},
\]
where $A_j = \{ x \in\m M\dvtx \| x - \tau\| = 2j \}, 0 \leq j \leq k$.
Since $k \leq s/2$, the expression in the above display can be bounded
above by $k {s \choose k} {p-1 \choose k}$. One thus has $|\m M| = {p-1
\choose s} \leq m k {s \choose k} {p-1 \choose k}$. Using $(n/r)^r \leq
{n \choose r} \leq(ne/r)^r$ for $0 \leq r \leq n/2$, we obtain $m \geq
\exp(C s \log p)$ for some constant $C > 0$. Also, clearly $m \leq|\m
M| \leq\exp(C_1 s \log p)$.
\end{pf*}

\begin{pf*}{Proof of Lemma \protect\ref{lem:shrinkage_conc}}
Let $\delta= \varepsilon/p$. To lower-bound $\bbP(\Vert\theta-
\theta _0\Vert_2 < \varepsilon)$, we first obtain a lower bound
conditioned on the
hyper parameters $\tau$ and~$\gamma$:
%
\renewcommand{\theequation}{\Alph{section}.\arabic{equation}}
\setcounter{equation}{0}
\begin{eqnarray}
\label{eq:prior_conc_cond} && \bbP\bigl(\Vert\theta- \theta_0\Vert_2 <
\varepsilon\mid\tau, \gamma\bigr)
\nonumber
\\
&&\qquad \geq\bbP\bigl(\vert\theta_j\vert \leq\delta\ \forall  j \in
S_0^c \mid \tau, \gamma\bigr)   \bbP\bigl(\Vert
\theta_{S_0} - \theta_{0S_0}\Vert_2 < \varepsilon/2
\mid\tau, \gamma\bigr)
\\
&&\qquad = \biggl[\prod_{j \in S_0^c} \bigl(1 - e^{-{\delta}/{\psi
_j}}
\bigr) \biggr] \times\bbP\bigl(\Vert\theta_{S_0} - \theta_{0S_0}
\Vert_2 < \varepsilon /2 \mid\tau, \gamma\bigr).\nonumber
\end{eqnarray}

Let $\tilde{\gamma} = (\gamma_1, \ldots, \gamma_{p-1})^{\T}$ and
$\gamma
_p = 1 - \sum_{j = 1}^{p-1} \gamma_j$. We now have to integrate out
$\tau$ and $\tilde{\gamma}$ in (\ref{eq:prior_conc_cond}). By a
relabeling of indices, we can always make sure that the $p$th index
lies in $S_0$. Let $S_1 = S_0 \setminus\{p\}$ so that $S_0^c \cup S_1
= \{1, \dots,\break  p-1\}$.
Fix numbers $a, b \in(0,1)$ with $b = 4a$. Observe that if $\tau\in
[2s, 4s], \gamma_j \tau\leq\frac{\delta}{\log(p/s)}\ \forall  j
\in S_0^c$ and $\gamma_j \tau\in[a, b]\ \forall  j \in S_1$, then
for $\varepsilon< b/2$,
%
\begin{equation}
\label{eq:simplex} \sum_{j=1}^{p-1}
\gamma_j = \sum_{j \in S_0^c}
\gamma_j + \sum_{j
\in
S_1}
\gamma_j \leq\varepsilon+ \frac{(s-1) b}{2s} \leq b < 1.
\end{equation}
Define $\mathcal{B} \subset\Delta_0^{p-1} \times\mathbb{R}^+$ such that
%
\begin{eqnarray}
\label{eq:B} &&\mathcal{B} = \biggl\{ (\gamma, \tau)\dvtx 0 \leq
\gamma_j \tau\leq \frac
{\delta}{\log(p/s)}\ \forall  j \in
S_0^c; \gamma_j \tau\in[a, b]\ \forall
j \in S_1,
\nonumber
\\[-8pt]
\\[-8pt]
\nonumber
&& \hspace*{230pt}\tau\in[2s, 4s] \biggr\}.
\end{eqnarray}
Clearly, $\mathcal{B}$ is a measurable subset of $\Delta_0^{p-1}
\times
\mathbb{R}^+$. For a fixed $\tau$ in the interval $[2s, 4s]$, the
section $\mathcal{A}_{\tau} \subset\Delta_0^{p-1}$ is given by
%
\begin{equation}
\label{eq:atau} \mathcal{A}_{\tau} = \biggl\{ 0 \leq\gamma_j
\leq\frac{\delta
}{\log
(p/s) \tau}\ \forall  j \in S_0^c;
\gamma_j \in \biggl[\frac
{a}{\tau
}, \frac{b}{\tau} \biggr]\
\forall  j \in S_1 \biggr\}.
\end{equation}

Thus,
%
\begin{eqnarray}\qquad
\label{eq:prior_conc_marg1} \bbP\bigl( \Vert\theta- \theta_0\Vert_2 <
\varepsilon\bigr) & = &\int_{ (\tau,
\tilde
{\gamma}) \in\mathbb{R}^+ \times\Delta_0^{p-1} } \bbP\bigl(\Vert \theta-
\theta_0\Vert_2 < \varepsilon\mid\tau, \tilde{\gamma}\bigr)
f_{\gamma}(d \tilde {\gamma}) f_{\tau}(d\tau)
\nonumber
\\[-8pt]
\\[-8pt]
\nonumber
& \geq&\int_{ (\tau, \gamma) \in\mathcal{B} } \bbP\bigl(\Vert\theta- \theta_0
\Vert_2 < \varepsilon\mid\tau, \gamma\bigr) f_{\gamma}(d \tilde {
\gamma }) f_{\tau}(d\tau).
\end{eqnarray}
We now substitute the lower bound for $\bbP(\Vert\theta- \theta
_0\Vert_2
< \varepsilon\mid\tau, \gamma)$ from (\ref{eq:prior_conc_cond}) in
(\ref
{eq:prior_conc_marg1}) and lower-bound the two terms on the right-hand
side of (\ref{eq:prior_conc_cond}) individually.

For the first term, observe that for $(\tau, \gamma) \in\mathcal{B}$,
$\prod_{j \in S_0^c} (1 - e^{-\delta/\psi_j}  ) \geq(1- s/p)^{p-s}$.

To tackle the second term, we make use of Lemma~\ref{lem:doubleconc}.
By definition, $\psi_j \in[a, b]$ for all $j \in S_1$ whenever $(\tau,
\gamma) \in\mathcal{B}$. Further, along the lines of \eqref
{eq:simplex}, $\sum_{j=1}^{p-1} \gamma_j \in[a/8, b]$, and hence
$\gamma_p \in[1 - b, 1- a/8]$ on $\mathcal{B}$. Hence, $\psi_p \in
[2s(1 - b), 4s(1- a/8)]$.
Since $a, b$ are constants, by a slight abuse of notation, we shall
assume $\psi_j \in[a, b]$ for all $j \in S_1$ and $\psi_p \in[2sa,
4sb]$ on $\mathcal{B}$.
It thus follows from Lemma~\ref{lem:doubleconc} that
\begin{eqnarray*}
&& \bbP\bigl(\bigl\Vert\Pi_{S_0}(\theta) - \Pi_{S_0}(
\theta_0)\bigr\Vert_2 < \varepsilon/2 \mid\tau, \gamma\bigr)
\\
&&\qquad \geq\exp \biggl\{-\frac{C_1}{a^2}\sum_{j \in S_0}
\vert\theta _{0j}\vert^2 - C_2s - s\bigl\vert\log
\bigl\{\varepsilon/(2b\sqrt{s})\bigr\}\bigr\vert \biggr\}.
\end{eqnarray*}
%
We conclude that for $(\tau, \gamma) \in\mathcal{B}$, the integrand in
(\ref{eq:prior_conc_marg1}) can be bounded below as follows:
%
\begin{eqnarray}
\label{eq:prior_conc_marg2}  &&\bbP\bigl(\Vert\theta- \theta_0\Vert_2 <
\varepsilon\mid\tau, \gamma\bigr)
\nonumber
\\[-8pt]
\\[-8pt]
\nonumber
&&\qquad \geq e^{-C s} \exp \biggl\{-
\frac{C_1}{a^2}\sum_{j \in S_0} \vert \theta
_{0j}\vert^2 -C_2s - s\bigl\vert\log\bigl\{
\varepsilon/(2b\sqrt{s})\bigr\}\bigr \vert \biggr\},
\end{eqnarray}
where the last inequality uses $(1 - x)^{1/x} \geq1/(2e)$ for $0 \leq
x \leq1/2$ and $C = \log(2e)$.
It thus remains to obtain a lower bound to
%
\begin{equation}
\label{eq:prior_Bset} \bbP(\mathcal{B}) = \int_{ (\tau, \gamma) \in\mathcal{B} }
f_{\gamma
}(d \tilde{\gamma}) f_{\tau}(d\tau) = \int
_{\tau= 2s}^{4s} \bbP( \mathcal{A}_{\tau} \mid
\tau) f_{\tau
}(d\tau).
\end{equation}
Now, since $\gamma\sim\operatorname{Dir}(\alpha/p, \ldots, \alpha/p)$,
recalling the definition of $\mathcal{A}_{\tau}$ from (\ref{eq:atau})
and using (\ref{eq:simplex}),
%
\begin{eqnarray}
\label{eq:pr_atau} &&\bbP(A_{\tau} \mid\tau) \nonumber\\
&&\qquad= \frac{\Gamma(\alpha)}{\Gamma(\alpha/p)^p} \int
_{\tilde{\gamma} \in\mathcal{A}_{\tau}} \Biggl[\prod_{j=1}^{p-1}
\gamma_j^{\alpha/p-1} \Biggr] \Biggl(1 - \sum
_{j=1}^{p-1} \gamma_j
\Biggr)^{\alpha/p - 1} \,d\gamma_1 \cdots d\gamma_{p-1}
\nonumber
\\[-8pt]
\\[-8pt]
\nonumber
&&\qquad \geq C_p (1 - b)^{\alpha/p-1} \int_{\tilde{\gamma} \in\mathcal
{A}_{\tau}}
\biggl[\prod_{j \in S_1} \gamma_j^{\alpha/p-1}
\biggr] \times \biggl[\prod_{j \in S_0^c}
\gamma_j^{\alpha/p-1} \biggr] \,d\gamma_1 \cdots d
\gamma_{p-1}
\\
&&\qquad\geq C_p (1 - b)^{\alpha/p-1} \biggl\{\frac{\delta}{\log
(p/s)} \biggr\}
^{\alpha(p-s)/p} \biggl\{ \biggl(\frac{b}{\tau} \biggr)^{\alpha/p} -
\biggl( \frac{a}{\tau} \biggr)^{\alpha/p} \biggr\}^{s-1},\nonumber
\end{eqnarray}
where
\begin{eqnarray}
\label{eq:cpineq}  C_p &=& \frac{\Gamma(\alpha)}{\Gamma(\alpha/p)^p} \biggl(\frac
{p}{\alpha
}
\biggr)^{p-1}
\nonumber
\\
& =& \exp\bigl\{ \log\Gamma(\alpha) + (p-1)\log(p/\alpha) - p \log \Gamma (
\alpha/p) \bigr\}
\nonumber
\\[-8pt]
\\[-8pt]
\nonumber
& \geq&\exp\bigl\{ \log\Gamma(\alpha) - \log\Gamma(\alpha/p) \bigr\}
\nonumber
\\
& \geq&\exp\bigl\{ \log\Gamma(\alpha) - \log(p/\alpha) \bigr\}\nonumber
\end{eqnarray}
with the last two inequalities using $\Gamma(x) \leq1/x$ for all $x
\in(0,1)$. Moreover, since $b \geq4a$, we have for $\tau\in[2s, 4s]$,
%
\begin{eqnarray}
\label{eq:contr_s0}  \biggl\{ \biggl(\frac{b}{\tau} \biggr)^{\alpha/p} -
\biggl( \frac
{a}{\tau
} \biggr)^{\alpha/p} \biggr\}^{s-1} &\geq&
\biggl\{ \biggl(\frac
{b}{4s} \biggr)^{\alpha/p} - \biggl(
\frac{a}{2s} \biggr)^{\alpha/p} \biggr\}^{s-1}
\nonumber\hspace*{-35pt}
\\[-4pt]
\\[-12pt]
\nonumber
& \geq& \biggl(\frac{b}{4s} \biggr)^{(s-1)\alpha/p} \biggl[ 1 - \exp \biggl
\{ - \frac{\alpha}{p} \log(2b/a) \biggr\} \biggr].\hspace*{-35pt}
\end{eqnarray}
Equations (\ref{eq:cpineq}) and (\ref{eq:contr_s0}), in conjunction
with the fact that $1 - e^{-x} \geq x/2$ for $x \in(0,1)$ implies that
the expression in
(\ref{eq:pr_atau}), and thus $\bbP(\mathcal{A}_{\tau} \mid\tau)$ in
(\ref{eq:prior_Bset}), is bounded below by
%
\begin{eqnarray}
\label{eq:pr_atau1}\qquad && \bbP\bigl(\mathcal{A}_{\tau} \mid\tau\bigr)
\nonumber
\\[-8pt]
\\[-8pt]
\nonumber
&&\qquad \geq C \exp \biggl\{
\frac
{\alpha
(p-s)}{p} \log\frac{\delta}{\log(p/s)} - \log\frac{p}{\alpha} -
\frac
{1}{\log(b/2a)} \log\frac{p}{\alpha} \biggr\}
\end{eqnarray}
for some constant $C > 0$. Finally, (\ref{eq:prior_conc_marg2}) and
(\ref{eq:pr_atau1}) substituted into (\ref{eq:prior_conc_marg1})
gives us
\[
\bbP\bigl(\Vert\theta- \theta_0\Vert_2 < \varepsilon\bigr) \geq
\bbP\bigl[\tau\in(2s, 4s) \bigr] e^{ - C \max\{\Vert\theta_0\Vert_2^2, s \log(s/\varepsilon),
\log p\} }.
\]
The proof of Lemma~\ref{lem:shrinkage_conc} is completed upon
observing that
$\bbP[\tau\in(2s, 4s)  ] \geq e^{-Cs}$.
\end{pf*}

%
%

\begin{pf*}{Proof of Lemma \protect\ref{lem:shrinkage_dim}}
Without loss of generality, we provide the proof for $\alpha= 1$.
Lemma IV.3 of \cite{zhou2012negative} implies that under (\ref{eq:shprior_prop}),
%
\begin{equation}
\label{eq:mingyuan} \theta_j \mid\psi_j \stackrel{
\mathrm{ind.}} \sim\operatorname{DE}(\psi_j),\qquad \psi_j
\stackrel{\mathrm{i.i.d.}}\sim\operatorname{Ga}(1/p, 1/2).
\end{equation}
By \eqref{eq:mingyuan}, $\theta_j$'s are independent and identically
distributed, so that\break  $\vert\operatorname{supp}_{\delta}(\theta)\vert \sim
\operatorname
{Binomial}(p, \zeta)$, with $\zeta:= \bbP(| \theta_1| > \delta)$. We
first show that $\zeta\lesssim\log p/p$ for $\delta= \varepsilon/p$.
Observe that
%
\begin{eqnarray}\label{eq:ig1}\quad
&&\bbP\bigl( \vert\theta_1\vert > \delta\bigr) \nonumber\\
&&\qquad= \frac{(1/2)^{1/p}}{\Gamma
(1/p)}\int
_{0}^{\infty} e^{-\delta/x} x^{1/p -1}
e^{-x/2} \,dx
\nonumber
\\[-8pt]
\\[-8pt]
\nonumber
&&\qquad= \frac{(1/2)^{1/p}}{\Gamma(1/p)} \biggl\{ \int_{0}^{4 \delta}
e^{-\delta/x} x^{1/p -1} e^{-x/2} \,dx + \int_{4 \delta}^{\infty}
e^{-\delta/x} x^{1/p -1} e^{-x/2} \,dx \biggr\}
\\
&&\qquad\leq \frac{(1/2)^{1/p}} {\Gamma(1/p)} \biggl\{C + \int
_{4 \delta
}^{\infty} \frac{e^{-x/2}}{x} \,dx \biggr\} \leq
\frac{(1/2)^{1/p}} {\Gamma(1/p)} \biggl\{ C + \int_{2
\delta}^{\infty}
\frac{e^{-t}}{t} \,dt \biggr\}.\nonumber  
\end{eqnarray}
Using a bound for the incomplete gamma function from Theorem~2 of \cite
{alzer1997some},
%
\begin{equation}
\label{eq:ig2} \int_{2\delta}^{\infty} \frac{e^{-t}}{t}
\,dt \leq-\log\bigl(1 - e^{-2
\delta
}\bigr) \leq-\log(\delta),
\end{equation}
for $\delta$ small. Since $\Gamma(1/p) \geq p/2$ for large $p$, and $C
+ \log(1/\delta) \leq2 \log(1/\delta)$ for $p$ large, we have
$\bbP(|
\theta_1| > \delta) \leq\log(1/\delta)/p \lesssim\log p/p$; the last
inequality follows since $\varepsilon> 1/p^B$ implies $\delta\geq1/p^{B+1}$.

A version of Chernoff's inequality for the binomial distribution \cite
{hagerup1990guided} states that for $B \sim\operatorname{Binomial}(p,\zeta)$
and $\zeta\leq a < 1$,
%
\begin{equation}
\label{eq:chernoff} \bbP(B > a p) \leq \biggl\{ \biggl(\frac{\zeta}{a}
\biggr)^a e^{a-\zeta} \biggr\}^p.
\end{equation}
In \eqref{eq:chernoff}, set $a = A s/p$. Since $s \gtrsim\log p$, we
can ensure $\zeta\leq a$ by choosing $A$ larger than some constant.
Hence, by \eqref{eq:chernoff}, $\bbP(\vert\operatorname{supp}_{\delta
}(\theta)\vert
> As) \leq(e \zeta/a)^{As} \leq e^{-A \log(A/eC) s}$.
\end{pf*}
\begin{pf*}{Proof of Lemma \protect\ref{lem:largedev_ps}}
Recall $\theta_j \mid\gamma, \tau\sim\operatorname{DE}(\gamma_j \tau)$ for
$1\leq j \leq p$. Let $X_j = \theta_j/(\gamma_j \tau)$, so that $X_j
\mid\gamma, \tau\sim\operatorname{DE}(1)$ independently. Let $\psi_j =
\gamma
_j \tau$ and fix $t > 1$. Using a Bernstein-type tail inequality for
subexponential random variables (Proposition~5.16 of \cite
{vershynin2010introduction}),
\begin{eqnarray*}
\bbP \Biggl(\sum_{j=1}^p |
\theta_j| > t \Bigm|\gamma, \tau \Biggr) &=& \bbP \Biggl(\sum
_{j=1}^p |\psi_j X_j |
> t \Bigm|\gamma, \tau \Biggr)\\
& \leq&\exp \biggl\{-C \min \biggl(\frac{t^2}{\Vert\psi\Vert_2^2}, \frac
{t}{\Vert\psi\Vert_{\infty}} \biggr)
\biggr\} \\
&\leq&\max \bigl\{e^{- C t^2/\tau^2}, e^{-C t/\tau} \bigr\}.
\end{eqnarray*}
The last inequality in the above display uses $\| \gamma\|_2^2 \leq\|
\gamma\|_1 = 1$ and $e^{-c/x}$ is increasing in $x$. Fix $t \geq1$.
Since $\tau\sim\operatorname{Exp}(1/2)$, $\bbP(\tau> \sqrt{t}) \leq e^{-C
\sqrt{t}}$. Also, $\max_{0 \leq x \leq\sqrt{t}} \max \{e^{- C
t^2/x^2}, e^{-C t/x} \} \leq e^{-C \sqrt{t}}$. The result follows by
noting that
\begin{eqnarray*}
\bbP\bigl(\| \theta\|_1 > t\bigr) &\leq&\int_{x = 0}^{\sqrt{t}}
\max \bigl\{ e^{- C
t^2/x^2}, e^{-C t/x} \bigr\} f_{\tau}(x) \,dx
+ \bbP(\tau> \sqrt{t}) \\
&\leq &2 e^{- C \sqrt{t}}.
\end{eqnarray*}
\upqed\end{pf*}
\end{appendix}

\section*{Acknowledgements}
The authors would like to thank two anonymous
referees, an Associate Editor and the Editor for their thoughtful
comments on previous versions of the paper which has helped improve our
exposition. We also thank Steven Finch for careful proofreading of an
initial draft of the paper.

%

%



\printaddresses

\end{document}